\documentclass[review,authoryear]{elsarticle}
\usepackage{amsfonts}
\usepackage{hyperref}
\usepackage{graphicx}
\usepackage{float}
\usepackage{amsmath, amssymb}
\usepackage{booktabs}
\usepackage{caption}
\usepackage{geometry}
\usepackage{longtable}
\usepackage{array}
\usepackage{multirow}  
\usepackage{xcolor}
\usepackage{mathtools}
\usepackage{amsthm}
\usepackage{natbib}
\usepackage[T1]{fontenc}
\usepackage[utf8]{inputenc}
\usepackage[polish,english]{babel}

\usepackage{todonotes}
\usepackage{easyReview}

\newtheorem{theorem}{Theorem}
\newtheorem{remark}{Remark}
\newtheorem{corollary}{Corollary}
\newtheorem{proposition}{Proposition}
\def\R{\mathbb{R}}
\newcommand{\RL}{\mathbb{R}}

\newcommand{\Exp}{\mathbb{E}}

\newcommand{\Prob}{\mathbb{P}}

\newcommand{\Var}{{\bf Var}}

\newcommand{\di}{\mathrm{d}}
\newcommand{\ind}{\mathbb{I}}

\def\zm#1{{\textcolor{red}{#1}}}

\def\zm#1{{#1}}

\DeclarePairedDelimiterXPP\pk[1]{\mathbb{P}}( ){}{ #1}

\def\akn#1{\begin{align} #1 \end{align}}

\def\bqn#1{\begin{eqnarray} #1 \end{eqnarray}}

\usepackage{comment}

\makeatletter
\def\ps@pprintTitle{%
 \let\@oddhead\@empty
 \let\@evenhead\@empty
 \let\@oddfoot\@empty
 \let\@evenfoot\@empty}
\makeatother

\begin{document}
\begin{frontmatter}
  
\title{Distributional and Extremal Behaviour of Brownian Motion \\
with Exponential Resetting}
 
\author[label1]{Krzysztof D\c{e}bicki}
\ead{Krzysztof.Debicki@math.uni.wroc.pl}

\author[label2]{Enkelejd  Hashorva}
\ead{Enkelejd.Hashorva@unil.ch}

\author[label3]{Zbigniew Michna\corref{cor1}}
\ead{zbigniew.michna@pwr.edu.pl}
\cortext[cor1]{Corresponding author.}

\affiliation[label1]{
organization={Mathematical Institute, University of Wroc\l aw},
          addressline={pl. Grunwaldzki 2/4},
          city={Wrocław},
          postcode={50-384},
          country={Poland}
}

\affiliation[label2]{
organization={Department of Actuarial Science, University of Lausanne},
addressline={Chamberonne}, 
postcode={1015},
city={Lausanne},
country={Switzerland}}

\affiliation[label3]{
organization={Department of Operations Research and Business Intelligence,
Wroc{\l}aw University of Science and Technology},
          addressline={Wybrzeże Stanisława Wyspiańskiego 27},
          city={Wrocław},
          postcode={50-370},
          country={Poland}}

\begin{abstract}
 
We study the distributional and asymptotic properties of the supremum of Brownian motion with drift and exponential resetting. We obtain an explicit renewal-type formula for the distribution of the supremum and then derive an approximation for its survival function. Moreover, we find the asymptotics of the tail distribution of the \zm{supremum} and infimum. We also consider the stationary case and give a new explicit expression for the \zm{two-dimensional distributions} of such processes.  
\end{abstract}
 

\end{frontmatter}

\section{Introduction}
Stochastic resetting is a mechanism that appears in many everyday situations. Consider the act of looking for something - such as trying to spot a familiar face in a large crowd or searching for your lost keys at home. When the search drags on without success, it often feels natural to return to the starting point and start again. Comparable behaviour is found elsewhere: a browser that reloads a frozen webpage, an algorithm that restarts after following an unproductive path, or a reader who returns to the beginning of a difficult paragraph to better grasp its meaning. In each case, the act of resetting interrupts aimless wandering and helps the process remain efficient.

From a mathematical standpoint, such behaviours can be described using Brownian motion with stochastic resetting. In a standard Brownian search, the average time required to locate a target - known as the mean first passage time - may grow without bound or become extremely large. By contrast, when the process is interrupted by resets occurring at a fixed rate, the system reaches a non-equilibrium steady state, and the mean first passage time typically becomes finite, sometimes even minimized. This illustrates why restart mechanisms, both in theory and in real-world settings, can improve efficiency - whether in physical search dynamics, computational procedures, or everyday decision-making.

Stochastic resetting appears in many physical phenomena. In the article \cite{Lenzi2022Transient}, the authors study diffusive dynamics in heterogeneous environments subject to stochastic resetting. They show that resetting profoundly affects transient anomalous diffusion, leading to crossovers between subdiffusive, normal, and superdiffusive regimes. Brownian motion with stochastic resetting emerges as a particular limiting case of the general framework, illustrating how resetting can suppress anomalous behaviour and influence the existence of stationary states.

\zm{A related Brownian model} 
with stochastic resetting to randomly distributed positions in the context of search and target-finding problems has been studied in \cite{ToledoMarin2023First}. This model is motivated by physical and biological processes in which a diffusing particle intermittently restarts its motion from new locations, such as intermittent search strategies or molecular transport. The work analyzes how the resetting position distribution influences first-passage-time statistics and information-theoretic measures, with resetting to a fixed point recovered as a special case.

In the paper of \cite{Evans2011DiffusionResetting}, \zm{a model in which a Brownian particle is intermittently reset to a fixed position $x_R$ at a fixed rate $\lambda>0$ is introduced and studied}. \zm{It is shown there} that this mechanism profoundly changes the long-time behaviour of the process: whereas pure diffusion spreads without bound, the addition of resetting produces a non-equilibrium stationary state whose distribution is Laplace-shaped, sharply peaked around the reset point. 
\zm{\cite{Evans2011DiffusionResetting} also investigate} first-passage properties, demonstrating that the presence of resetting makes the mean first-passage time to a target finite, in contrast to the divergent mean time in ordinary diffusion. Moreover, an explicit expression for the Laplace transform of the distribution of the first-passage time and the optimal resetting rate that minimises the expected time to reach the target is analyzed. 
However, the survival probability of the first-passage time is implicitly given by the inverse Laplace transform (see also \cite{Hartmann2024}).

In \cite{EvansMajumdar2011JPA}, several generalisations of their previous model are considered, including: (i) a space-dependent
resetting rate; (ii) resetting to a random position drawn from a resetting
distribution; and (iii) a spatial distribution for the absorbing target.
We refer to, e.g.,  \cite{EvansMajumdarSchehr2020} for the analysis of different types of resetting protocols, such as Poissonian, non-Poissonian, and memory-dependent resetting.

\zm{In \cite{GuptaPlataPal2020JarzynskiReset}, the authors study an overdamped Brownian particle in a time-dependent potential subject to stochastic resetting and analyze the resulting work fluctuations under nonequilibrium driving. They show that resetting leads to a universal Gaussian form of the work distribution at long times and investigate the validity of the Jarzynski equality, which is generally violated under the finite-time resetting unless specific conditions on the driving protocol are satisfied.
}

The article \cite{MonteroMasoPuigdellosasVillarroel2017} investigates the effects of reset events on continuous-time random walks (CTRWs). It provides a historical overview, introduces a general analytical framework, and presents new results for a monotonic CTRW with constant drift that may change direction after each reset. The authors derive the transition probability density for arbitrary jump distributions and drift speeds, demonstrate the emergence of a stationary distribution, and provide formulas for the mean first-passage time, highlighting conditions under which it can be minimised. Analytical results are validated through Monte Carlo simulations, showing excellent agreement with theoretical predictions.

\zm{Other types of resetting appear in \cite{ZhangHadjiliadis2012DrawdownsCrash}. Their motivation comes from the financial risk management for modelling and quantifying extreme portfolio losses, as well as from the pricing and hedging of insurance-type derivatives that protect against rapid market crashes.}

\zm{In \cite{Silvestrov2021DampingMC}, a resetting-like mechanism is used in perturbed Markov chains to model systems with random reinitialization or global mixing, where the dynamics occasionally jump according to an auxiliary distribution, which helps ensure convergence to a well-defined stationary regime.}

In \cite{MagdziarzTazbierski2022}, a general stochastic representation for processes with resetting is introduced, modelling them as jump-diffusion systems that intermittently restart from a prescribed point. This framework provides both analytical tools and Monte Carlo simulation methods for studying resetting dynamics. As a benchmark, the authors derive fundamental properties of Brownian motion with Poissonian resetting, including the Itô formula, the probability density, the moment-generating and characteristic functions, moments of all orders, the Fokker–Planck equations, and the infinitesimal generator and its adjoint. The results are further extended to time-inhomogeneous Poissonian resetting, yielding a general framework for stochastic processes with random resetting.

The recent article \cite{TazbierskiMetzlerMagdziarz2025} 
develops a series-representation framework for stochastic processes subject to renewal resetting. It treats both complete resetting, which erases the system’s memory, and incomplete resetting, which only shifts the process position. Using this approach, the authors derive explicit expressions for the joint two-time probability density and the autocorrelation function, which were previously unavailable. The framework is further applied to first-passage problems, yielding general results for complete resetting. The theory is applied to the Brownian motion and the scaled Brownian motion under Poissonian and mixture-exponential resetting, with the latter also serving as an approximation scheme for more general resetting time distributions.

In this paper, we derive the two-dimensional joint distribution of Brownian motion with drift and exponential resetting, and its limiting two-dimensional distribution as time tends to infinity (corresponding to the stationary regime). 
A substantial part of this contribution is dedicated to the extremal properties of Brownian motion with resetting.

\zm{Section \ref{sec2} is devoted to the distributional properties of Brownian motion with resetting. Specifically, two-dimensional joint cumulative distribution function and their limits are found.}

In Section \ref{sec3}, we derive the exact formula for the distribution of the supremum, its tail distribution asymptotics,
and quantities straightforwardly connected with the notion of the first-passage time.
We note that so far, only the Laplace transform of the first-passage-time distribution has been known in explicit form (see \cite{Evans2011DiffusionResetting}).  
Our findings are complemented by the asymptotic analysis of the tail distribution of the finite time infimum and the value at the last point.

Section \ref{sec4} is dedicated to the analysis of supremum of the stationary Brownian motion with exponential resetting
over a given time interval, including 
its distribution function, exact asymptotics, and asymptotic properties of the joint distribution of the finite time supremum and the value at the last point.
\def\R{\mathbb{R}}

Additionally, in Section \ref{sec5} we present numerical examples that illustrate the derived theoretical results.
The proofs of all the theoretical results are postponed to the Appendix \ref{app}.

\section{Distributional properties of Brownian motion with resetting}\label{sec2} 
\begin{figure}[H]
\centering
\includegraphics[width=0.6\textwidth]{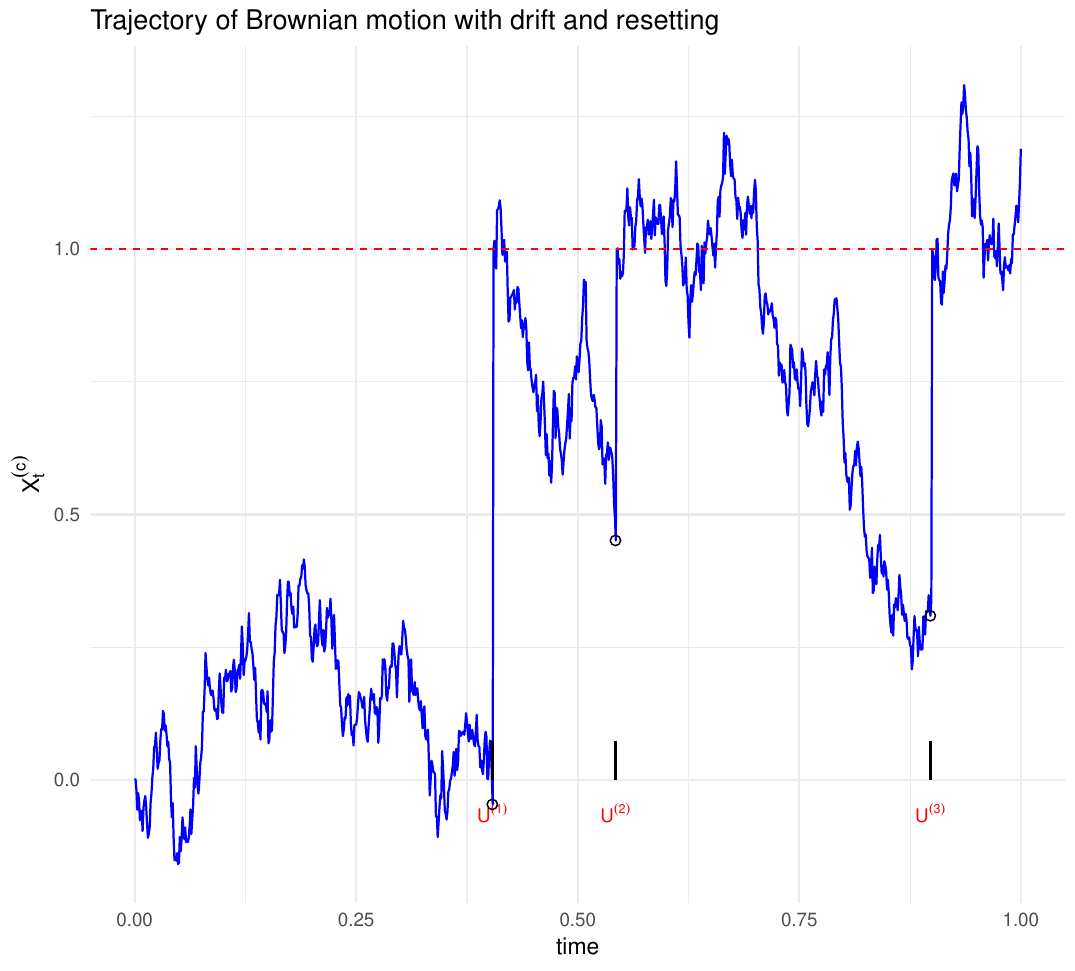}
    \caption{A trajectory of Brownian motion with resetting with $c=1$, $x_0=0$, $x_R=1$, $\lambda=2$. $U^{(1)},\, U^{(2)},\,U^{(3)}$ are epochs of resetting.}
\label{Pic:1}
\end{figure}

Let $(\Omega,\mathcal F,\Prob)$ be a probability space supporting an independent pair
$(N,W)$, where $(W_t)_{t\ge0}$ is a Brownian motion with $\Var(W_t)=\sigma^2 t$
and $(N_t)_{t\ge0}$ is a Poisson process with intensity $\lambda>0$,
which is independent of $W$.
\zm{To simplify notation, without loss of generality, in the rest of the paper we assume that $\sigma=1$.}
Let $(\mathcal F_t)_{t\ge0}$ be the usual augmentation of the natural filtration
generated by $(W,N)$.
Fix $c\in\R$, a reset level $x_R\in\R$, and an initial value $x_0\in\R$. If $x_0$ is random, we assume that it is  $\mathcal F_0$ measurable and independent of $(W,N)$. 
Define the drifted Brownian motion $W^{(c)}_t\coloneqq W_t-ct$.\\
We define Brownian motion with drift $-c$ and exponential resetting as the unique
càdlàg $(\mathcal F_t)$--adapted process
$X^{(c)}=(X^{(c)}_t)_{t\ge0}$ with $X^{(c)}_0=x_0$ satisfying the integral equation
\begin{equation}\label{eq:def-reset-sde}
X^{(c)}_t
=
x_0 + W^{(c)}_t
+ \int_{(0,t]}(x_R - X^{(c)}_{s-})\,dN_s,
\qquad t\ge0.
\end{equation}
\zm{For a trajectory of Brownian motion with exponential resetting see Fig. \ref{Pic:1}}. 
The above jump-diffusion representation for $c=0$ was introduced in  \cite{MagdziarzTazbierski2022}.  
Equivalently, $X^{(c)}$ can be defined pathwise as the càdlàg
$(\mathcal F_t)$-- adapted process
\begin{equation}\label{eq:tau-repr-x0}
X^{(c)}_t
\coloneqq x_R + \bigl(W^{(c)}_t - W^{(c)}_{\tau(t)}\bigr)
   + (x_0-x_R)\ind\{\tau(t)=0\},
\qquad t\ge0,
\end{equation}
where
\[
\tau(t)\coloneqq\sup\{s\in(0,t]: \Delta N_s=1\},\qquad t\ge0,
\]
with the convention $\tau(t)=0$ if $N_t=0$, and
where $\Delta N_s:=N_s-N_{s-}$ denotes the jump of $N$ at time $s$.\\
The equivalence between \eqref{eq:def-reset-sde} and
\eqref{eq:tau-repr-x0} follows by direct inspection of the dynamics
between successive jump times of $N$ and the reset rule at jump epochs.\\
As an immediate consequence of the pathwise construction {(\ref{eq:tau-repr-x0})}, the process
$X^{(c)}$ is a strong Markov process with respect to the filtration
$(\mathcal F_t)_{t\ge0}$, see \cite{AvrachenkovPiunovskiyZhang2013}.
In the rest of this section, we collect some basic distributional properties of $X^{(c)}_t$,
beginning with the distribution function of  $X^{(c)}_t$. We note that for $c=0$, $x_0=0$, the explicit form of the density function 
of $X^{(0)}_t$ can be found in \cite{MagdziarzTazbierski2022}; see also
the master equation \zm{(extension of Fokker-Planck equation)} in \cite{Evans2011DiffusionResetting}.

Using that
    $$
    \Prob(X^{(c)}_t\leq u)=e^{-\lambda t}\Prob(X^{(c)}_t\leq u | N_t=0)+(1-e^{-\lambda t})\Prob(X^{(c)}_t\leq u | N_t\geq 1)
    $$
and the fact that the time since the last jump
$t-U^{(N_t)}$ under the condition $N_t\geq 1$, 
where $U^{(1)}, U^{(2)},\ldots$ are the epochs of the consecutive jumps of the Poisson process, 
is a truncated exponential  
with density
$g_t(x)=\frac{\lambda e^{-\lambda x}}{1-e^{-\lambda t}}$ for $x\in (0,t)$, 
we immediately arrive at
\begin{eqnarray}
  \Prob(X^{(c)}_t\leq u)=e^{-\lambda t}\Prob(W_t-ct\leq u - x_0)+\lambda\int_0^t
    \Prob(W_s-cs\leq u-x_R)e^{-\lambda s}\di s\,.\label{oned}
\end{eqnarray} 
 
Importantly for the rest of this paper, the following convergence in distribution,  as $t\to \infty $, holds 
$$
X^{(c)}_t \rightarrow X^{(c)}_\infty \overset{d}{=}x_R+ \sqrt{S} W_1-cS ,   
$$
where $S$ is an exponential random variable with parameter $\lambda>0$ independent of the Wiener  process $W$.
The random variable $X^{(c)}_{\infty}$ has the asymmetric Laplace distribution (double exponential) whose density function 
has the following form
\begin{equation}\label{finfty}
f_{X^{(c)}_{\infty}}(x)=\frac{\lambda}{\sqrt{c^2+2\lambda}}
\exp\!\left(
\zm{-}c\,(x-x_R)
-\sqrt{c^2+2\lambda}\,
\lvert x-x_R\rvert
\right),
\qquad x\in\RL\,.
\end{equation}
Hence, for
$
\alpha
\coloneqq
\sqrt{c^2+2\lambda}\zm{+}c
>0$
and $u\geq x_R$
\begin{equation}\label{eq:tail-Y0}
\Prob\!\big(X^{(c)}_\infty>u\big)
=
\frac{\sqrt{c^2+2\lambda}\zm{-}c}
     {2\sqrt{c^2+2\lambda}} e^{-\alpha (u-x_R)}\,.
\end{equation}
Moreover,  we have 
\begin{equation}\label{inftydistr}
    \Prob(X^{(c)}_\infty\leq u)=\Prob(W_S^{(c)}\leq u-x_R)=\lambda\int_0^\infty
    \Prob(W_s-cs\leq u-x_R)e^{-\lambda s}\di s;
\end{equation}
compare also with eq. (2) in \cite{Evans2011DiffusionResetting}.
 
The joint probability density function of 
Brownian motion with a general renewal resetting was derived in \cite{TazbierskiMetzlerMagdziarz2025}. 
In the next theorem, 
we give a formula for the joint cumulative distribution function of the Brownian motion with exponential resetting,
and random $x_0$.

\begin{theorem}\label{jointdist_drift}
Assume that $X^{(c)}_0=x_0$, where $x_0$ is a real-valued random variable
independent of $(W,N)$.
For $0\le s<t$ and $u,w\in\RL$
\begin{align}
\Prob(X^{(c)}_s\le u, X^{(c)}_t\le w)
={}& e^{-\lambda t}\,
\Prob\big(W^{(c)}_s\le u-x_0,\, W^{(c)}_t\le w-x_0\big) \label{cdfj_drift}\\
&+ e^{-\lambda (t-s)}\lambda
\int_0^s
\Prob\big(W^{(c)}_{x}\le u-x_R,\, W^{(c)}_{t-s+x}\le w-x_R\big)
e^{-\lambda x}\,\di x \nonumber\\
&+\Bigg[
e^{-\lambda s}\Prob\big(W^{(c)}_s\le u-x_0\big)
+\lambda\int_0^s
\Prob\big(W^{(c)}_x\le u-x_R\big)
e^{-\lambda x}\,\di x
\Bigg]\nonumber\\
&\qquad\cdot\lambda\int_0^{t-s}
\Prob\big(W^{(c)}_y\le w-x_R\big)
e^{-\lambda y}\,\di y .\nonumber
\end{align}
\end{theorem}
For other distributional characteristics of $X^{(c)}_t$ such as the
explicit form of the one-dimensional probability density function, \zm{the two-dimensional density function}, the moment-generating function, the characteristic function and moments of all orders, see \cite{MagdziarzTazbierski2022}.

Taking $t=s+\delta$
and letting $s\to\infty$ in \eqref{cdfj_drift}, we get the following corollary. We note that 
the effect of the initial value $x_0$ vanishes asymptotically due to
the regenerative structure induced by resetting, and the limit is independent of $x_0$, even if we allow it to be random and independent of $(W,N)$.
\begin{corollary}\label{cor3}
For $\delta>0$ and $u,w\in \RL$,
\bqn{
\lim_{s\rightarrow\infty}\Prob(X^{(c)}_s\leq u, X^{(c)}_{s+\delta}\leq w)
&=&
e^{-\lambda \delta}
\Prob(W^{(c)}_{S}\leq u-x_R,\,
      W^{(c)}_{\delta+S}\leq w-x_R)
\label{eq:stat-joint-drift}\\
&+&
\Prob(W^{(c)}_S\leq u-x_R)
\Prob(W^{(c)}_S\leq w-x_R,\; S\leq \delta )
\nonumber\\
&=&
e^{-\lambda \delta}
\Prob(W^{(c)}_{S}\leq u-x_R,\,
      W^{(c)}_{\delta+S}\leq w-x_R)
\nonumber\\
&+&
\Prob(X^{(c)}_\infty\leq u)\,
\Prob(\sqrt{S}W_1-cS\leq w-x_R,\; S\leq \delta )\,.
\nonumber}
\end{corollary}

\def\supT{\sup_{t \in [0,1]}}

\section{Extremal behaviour of Brownian motion with resetting}\label{sec3}
In this section, we discuss both the distribution of the supremum and infimum of a Brownian motion with resetting. 
We first consider a distributional result for  $\sup_{t\in [0,T]} X^{(c)}_t$.
A natural motivation for the analysis of the supremum functional stems
from its connection with {\it the search problems} used in the stochastic modelling of such 
processes as, e.g., animals looking for food, proteins locating specific binding sites on DNA, chemical molecules diffusing and seeking reaction partners, debugging procedures identifying bugs in computer programs, randomized algorithms exploring a high-dimensional landscape to find a global minimum (see e.g. \cite{Viswanathan2011}, \cite{ToledoMarin2023First} and \zm{\cite{SinghPal2021ResetEVS}}). 
More precisely, let $u\in \R$ denote the location of a fixed target, and define the time $\tau^{(c)}$ to reach the target by
\begin{equation}\label{tau}
    \tau^{(c)}= \inf\{t>0: X^{(c)}_t>u\}\,.
\end{equation}
Then, the following relation holds
\begin{eqnarray}
\Prob(\tau^{(c)}>T)=\Prob( \sup_{t\in [0,T]}X^{(c)}_t\le u)\,. \label{equiv1}
\end{eqnarray}
 
The backward master equation for the probability of the first passage of zero for the Brownian motion (without drift) with resetting $x_R$ starting from any point $x_0$ is provided
in \cite{Evans2011DiffusionResetting} eq. (4), and it could be easily generalised to any hitting point. Moreover, the Laplace transform of this probability is given in an explicit form (see \cite{Evans2011DiffusionResetting} eq. (6), \zm{see also \cite{Hartmann2024}}, the inverse Laplace transform of it exists only in an implicit form \zm{of the Bromwich integral}).

Recall that for $c\in \RL$, $u\geq 0$, $T>0$ we set 
\begin{equation}\label{probF}
F(u,T)\coloneqq \Prob\left(\sup_{t \in [0,T] } (W_t -ct)\le u\right)=\Phi\left(\frac{u+cT}{\sqrt{T}}\right)+
e^{-2uc}\Phi\left(\frac{u-cT}{\sqrt{T}}\right)-e^{-2uc}\, ,
\end{equation}
where $\Phi$ is the  distribution function of the standard normal distribution. Below we set $\varphi=\Phi', \Psi=1-\Phi$. 
Note that $F(u,T)=0$ for all $u\le0$.
 
\begin{theorem}\label{th2}
For $u\in\R$, $c\in\R$, $x_0\in\R$ and $T>0$ we have
\akn{
\Prob\left(\sup_{t \in [0,T]} X^{(c)}_t \le u\right) 
&= e^{-\lambda T} F(u-x_0,T) \nonumber\\
&\quad + e^{-\lambda T} \sum_{n=1}^\infty \lambda^n 
\!\!\!\!\int\limits_{\substack{s_1 + s_2 + \ldots + s_n < T \\ s_i > 0}} 
\!\!\!\! F(u-x_0,s_1) 
\prod_{k=2}^n F(u-x_R, s_k) \nonumber \\
&\qquad\qquad \times F\big(u-x_R,\,T-s_1-\ldots-s_n\big)
\, \di s_1 \di s_2 \ldots \di s_n\,,
\label{mr2}
}
where we use the convention $\prod_{k=2}^1=1$.
\end{theorem}

\begin{remark} 

We note that for $x_0\leq u\leq x_R$ we have $\Prob\left(\sup_{t\in [0,T]}  X^{(c)}_t \leq u\right) = e^{-\lambda T}F(u-x_0,T)$.  
The integral in (\ref{mr2}) can be written as an iterated integral
$$
\int_0^T\di s_1\int_0^{T-s_1}\di s_2\ldots\int_0^{T-s_1-\ldots-s_{n-1}}\di s_n
\,F(u-x_0,s_1) 
\prod_{k=2}^n F(u-x_R, s_k) \nonumber F(u-x_R, T-s_1-\ldots-s_n)\,.
$$
\end{remark}


\begin{remark}
Using (\ref{mr2}) we obtain  
\begin{equation}\label{ineqF}
   e^{-\lambda T}[1-F(u-x_0,T)]\leq \Prob(\sup_{t\in [0, T]} X^{(c)}_t > u)\leq
   1-F(u-x_0,T)e^{-\lambda T[1-F(u-x_R, T)]}\,,
\end{equation}
where $F(u,T)$ is defined in \eqref{probF}.  
\end{remark}

Next, we shall investigate the asymptotic behavior of  
$\Prob(\sup_{t\in [0, T]} X^{(c)}_t > u)$ as $u\rightarrow\infty$. 
\zm{In the next result}, we suppose that $x_0=0$.  
Below we write $f(u)\sim g(u)$ as $u\rightarrow\infty$ for positive functions $f$ and $g$ if $\lim_{u\rightarrow\infty}\frac{f(u)}{g(u)}=1$.
\begin{theorem}\label{th.sup.as}
\begin{itemize}
    \item[i)] If   $x_R\le0$, then as $u\to\infty$ uniformly for  $T\in [a,b]$, with $0<a<b$
    $$\Prob(\sup_{t\in [0, T]} X^{(c)}_t > u)\sim 2e^{-\lambda T}\Psi\left(\frac{u+cT}{\sqrt{T}}\right). $$ 
    \item[ii)] If  $x_R>0$, then as $u\to\infty$
       $$\Prob(\sup_{t\in [0, T]} X^{(c)}_t > u)\sim 4\lambda T^2 e^{-\lambda T}u^{-2}\Psi\left(\frac{u+cT-x_R}{\sqrt{T}}\right). $$
\end{itemize}    
\end{theorem}

\zm{The infimum distribution of stochastic processes with resetting plays an important role in search processes, such as animal foraging, biomolecular search, and chemical kinetics. This connection is discussed by \cite{Majumdar2020EVS} (Chapter 4.2.5), who relate survival probabilities under resetting to the extreme value statistics of process minima}.
\zm{Thus} we consider next the behaviour of the infimum on $[T,T+ \Delta]$ when it stays above a given level $u$,
together with the additional constraint $X^{(c)}_{T}>v$.
If $u>x_R$, then at any reset time $\tau$ we have $X^{(c)}_\tau=x_R<u$,
implying for $T \in [0,\infty), \Delta>0$
\[
\Big\{\inf_{t\in[T,T+\Delta]}X^{(c)}_t>u,\; X^{(c)}_{T}>v\Big\}
\subseteq
\{N_{T+\Delta}-N_{T}=0\}.
\]
Conditioning on the last reset time prior to $T$ and exploiting the explicit
structure of Brownian motion with drift between successive reset times,
we obtain   
\begin{equation}\label{eq:inf-window-exact-v}
\begin{aligned}
\pk*{\inf_{t\in[T,T+\Delta]}X^{(c)}_t>u,\; X^{(c)}_{T}>v}
&= 
e^{-\lambda (T+\Delta)}\,
F_{T+\Delta}^{(v-x_0+cT)}(u-x_0,\Delta)
\\
&\quad
+
\int_{\Delta}^{T+\Delta}\lambda e^{-\lambda s}
F_{s}^{(v-x_R+c(s-\Delta))}(u-x_R,\Delta)\,ds\\
&=
e^{-\lambda (T+\Delta)}\,
F_{T+\Delta}^{(v-x_0+cT)}(u-x_0,\Delta)
\\
&\quad
+
e^{-\lambda \Delta}\int_{0}^{T}\lambda e^{-\lambda s}
F_{s+\Delta}^{(v-x_R+cs)}(u-x_R,\Delta)\,ds,
\qquad u>x_R,
\end{aligned}
\end{equation}
where 
$$
F_s^{(w)}(u,\delta)
\coloneqq
\pk*{\inf_{t\in[s-\delta,s]}(W_t-ct)>u,\; W_{s-\delta}>w},
\qquad u,w\in\R.\\
$$

Asymptotic results for the tail of the infimum of light-tailed random processes are in general difficult to obtain, see \cite{adler2014existence,chakrabarty2018asymptotic,DebickiHashorvaNovikov2026}
for more details.

\begin{theorem}\label{thm:BM-window-reset-general-x0}
Fix $ T>0, \Delta>0$  and set 
\[
K_{c,\Delta}
\coloneqq
\frac{2}{\sqrt{\Delta}}
\varphi\big(c\sqrt{\Delta}\big)
-2 c \Psi\big(  c  \sqrt{\Delta}\big)>0.
\]
Let $r\ge0$ and define  
\[
v=u+\frac{r}{u+cT},
\qquad
L(y)=
\begin{cases}
1, & y\le 0,\\[1mm]
e^{-y}(1+y), & y>0.
\end{cases}
\] 
\medskip
\noindent (i) If $x_0<x_R$, then 
as $u\to\infty$
\[
\pk*{\inf_{t\in[T,T+\Delta]}X_t^{(c)} >u,\; X_{T}^{(c)} >v}
\sim
2\lambda\,e^{-\lambda (T+\Delta)}K_{c,\Delta}\,L(r/T)
\frac{T^3}{u^3}
\Psi\Big(\frac{u-x_R+cT}{\sqrt{T}}\Big).
\]

\medskip
\noindent (ii) If $x_0\ge  x_R$, then as $u\to\infty$
\[
\pk*{\inf_{t\in[T,T+\Delta]}X_t^{(c)} >u,\; X_{T}^{(c)} >v}
\sim
e^{-\lambda (T+\Delta)}\,
K_{c,\Delta}
L(r/T)\frac{T}{u}
\Psi\!\Big(\frac{u-x_0+cT}{\sqrt{T}}\Big).
\]
\end{theorem}

\def\eps{\varepsilon}

\section{Stationary Brownian motion with resetting}\label{sec4}

In this section, we analyse a
Brownian motion with exponential times of resetting, where $x_0$ is random. Specifically, we consider 
$x_0=X^{(c)}_\infty$ with the Laplace distribution as in eq. (\ref{finfty}) \zm{and independent of $(W,N)$.} Unlike the previous sections, we denote this process by $Y^{(c)}_t, t\ge 0.$ In \cite{AvrachenkovPiunovskiyZhang2013} basic properties of Markov processes with Poissonian resetting and a general form of the invariant measure are derived.  
 However, the article does not discuss stationary Markov processes.
 
\begin{proposition}\label{prop.stat}
     $Y^{(c)}_t,\,t\geq 0$ is stationary with 
    the joint distribution function of $(Y_0^{(c)},Y_\delta^{(c)}),\delta>0 $ 
    given in \eqref{eq:stat-joint-drift}.
\end{proposition}
\newcommand{\equaldis}{\stackrel{d}{=}}
As a direct conclusion from Theorem \ref{th2}, 
for $c,u\in\R$ and $T>0$, we get 
\akn{
\Prob\left(\sup_{t \in [0,T]} Y^{(c)}_t \le u\right) 
&= e^{-\lambda T} \Exp F(u-X^{(c)}_\infty,T) \nonumber\\
&\quad + e^{-\lambda T} \sum_{n=1}^\infty \lambda^n 
\!\!\!\!\int\limits_{\substack{s_1 + s_2 + \ldots + s_n < T \\ s_i > 0}} 
\!\!\!\! \Exp F(u-X^{(c)}_\infty,s_1) 
\prod_{k=2}^n F(u-x_R, s_k) \nonumber \\
&\qquad\qquad \times F\big(u-x_R,\,T-s_1-\ldots-s_n\big)
\, \di s_1 \di s_2 \ldots \di s_n\,,\nonumber
}
where $F$ is defined in (\ref{probF}), $X^{(c)}_\infty$ has distribution given in (\ref{finfty}) and we use the convention $\prod_{k=2}^1=1$.

Next, we derive an asymptotic behaviour of the supremum distribution for the process $Y_t\coloneqq Y_t^{(0)}$. 
\begin{theorem}\label{thasimpBMS} We have with   $\alpha=\sqrt{2\lambda}$  and $T>0$ as $u\rightarrow\infty$
    \begin{equation}\label{asimpBMS}
    \Prob(\sup_{t\in [0,T]}Y_t>u)\sim  e^{\alpha x_R}[\Phi(\alpha\sqrt{T})+ \lambda\int_0^T\Phi(\alpha\sqrt{s})\di s]\exp(-\alpha u).
    \end{equation} 
\end{theorem}

In the following theorem, we shall derive the asymptotics of the joint distribution of the supremum and the last 
value of a stationary Brownian motion with resetting. 

\begin{theorem}\label{thasympSBMl}
As $u\rightarrow\infty$
$$
\Prob(\sup_{t\in [0,T]}{Y_t>u}, Y_T>uz)\sim\left\{    
    \begin{array}{ll}
    e^{\alpha x_R}[\Phi(\alpha \sqrt{T})+ \lambda\int_0^T\Phi(\alpha\sqrt{s})\di s]\exp(-\alpha u)& \mbox{if}\,\,\,z<0  \\
 e^{\alpha x_R}\Phi(\alpha \sqrt{T})\exp(-\alpha u)  &\mbox{if} \,\,\,0<z<1\,.
\end{array}
\right.
$$
\end{theorem}
\begin{remark}
Theorem \ref{thasympSBMl} does not cover the case  $z=0$. We conjecture that  
    $$
\Prob(\sup_{t\in [0,T]}Y_t>u, Y_T> 0)\sim
C_T\exp(-\alpha u)
    $$
    as $u\rightarrow\infty$,
where $e^{\alpha x_R}\Phi(\alpha\sqrt{T})\leq C_T\leq e^{\alpha x_R}[\Phi(\alpha\sqrt{T})+ \lambda\int_0^T\Phi(\alpha\sqrt{s})\di s]$.

Note that if $z\geq 1$, then   $\Prob(\sup_{t\in [0,T]}Y_t>u, Y_T>uz)=\Prob(Y_T>uz)=\frac{1}{2}e^{-\alpha(uz-x_R)}$ for $uz\geq x_R$.
\end{remark}

\section{Numerical examples}\label{sec5}
We shall begin with the numerical simulation of  the expectation of the first passage time $\tau\coloneqq \tau^{(0)}$ of the Brownian motion with resetting $X_t\coloneqq X^{(0)}_t$ and $x_0=0$ (see (\ref{tau})).
We note that \zm{for $x_R=0$ and $u>0$}
\begin{equation}\label{exfor}
\Exp\tau = \frac{1}{\lambda}\left(\exp\left(u\sqrt{2\lambda}\right)-1\right)
\end{equation}
and it attains its minimum for $\lambda^*\approx(1.59362)^2/(2u^2)$, see \cite{Evans2011DiffusionResetting}.

Combining Theorem \ref{th2}  with (\ref{equiv1}), for $u\in\R$
we arrive at 
\begin{eqnarray}
\Prob(\tau>T)&= &e^{-\lambda T}F(u,T)+
e^{-\lambda T} \sum_{n=1}^\infty \lambda^n 
\!\!\!\!\int\limits_{\substack{s_1 + s_2 + \ldots + s_n < T \\ s_i > 0}} 
\!\!\!\!F(u,s_1) 
\prod_{k=2}^n F(u-x_R, s_k)\nonumber \\
&&\qquad\qquad\qquad\times F(u-x_R, T-s_1-\ldots-s_n)
\, \di s_1 \di s_2 \ldots \di s_n\,.\label{pt}
\end{eqnarray} 
To compute integrals on the simplex given in (\ref{pt}) we shall use the following equality
$$
\int\limits_{\substack{s_1 + s_2 + \ldots + s_n < T \\ s_i > 0}} 
g(s_1, s_2,\ldots,s_n) 
\, \di s_1 \di s_2 \ldots \di s_n=\frac{T^n}{n!}\Exp[g(D_1, D_2,\ldots,D_n)]\,,
$$
where $(D_1, D_2,\ldots,D_n)$ is a random vector with uniform distribution on the simplex $s_1 + s_2 + \ldots + s_n < T,\,\, s_i > 0$, $i=1,2,\ldots,n$, $T>0$ (Dirichlet distribution, see e.g. \cite{johnson2000continuous}). If $(S_1, S_2,\ldots, S_{n+1})$ are independent unit exponential random variables then
$$
(D_1, D_2,\ldots,D_n)\stackrel{d}{=}T\left(\frac{S_1}{\sum_{k=1}^{n+1}S_k}, \frac{S_2}{\sum_{k=1}^{n+1}S_k},\ldots, \frac{S_n}{\sum_{k=1}^{n+1}S_k} \right)\,,
$$
see e.g. \cite{johnson2000continuous}. Thus, applying Monte Carlo method to compute the above integral, we get
$$
\int\limits_{\substack{s_1 + s_2 + \ldots + s_n < T \\ s_i > 0}} 
g(s_1, s_2,\ldots,s_n) 
\, \di s_1 \di s_2 \ldots \di s_n=\frac{T^n}{n!}\lim_{M\rightarrow\infty}
\frac{1}{M}\sum_{k=1}^M g(D^{(k)}_1, D^{(k)}_2,\ldots,D^{(k)}_n)\,,
$$
where $M$ is the number of realizations of the vector $(D_1, D_2,\ldots,D_n)$.

Next using formula (\ref{pt}) we compute $\Exp \tau$ for different values of the intensity $\lambda$, see Tab. \ref{Tab:1} ($\Exp \tau$ APPR). 
The integral $\Exp \tau= \int_0^\infty \Prob(\tau>s)\di s$ is approximated by
$\int_0^{T_{max}} \Prob(\tau>s)\di s$ with $T_{max}=30$ and the discretization step $e=0.01$. We take $n=60$ summands from the formula (\ref{pt}) and the integrals on the simplex are estimated using Monte Carlo method with $M=5000$ samples. The exact value of  $\Exp \tau$ is also calculated
using formula (\ref{exfor})  ($\Exp \tau$ EXACT in Tab.\ref{Tab:1}).
\begin{table}[H]
\centering
\resizebox{\textwidth}{!}{%
\begin{tabular}{l|l|l|l|l|l|l|l|l|l|l|l}
\hline
$\lambda$ & $0.1$& $0.269812
$& $0.769812$& $ 1.069812$&$1.169812$ & $1.269812^*$&$ 1.369812$ & $ 1.469812$&$1.769812$ &$2.269812$ &$4.269812$          \\
\hline
$\Exp\tau\,\, \mbox{APPR}$ & $5.315687$&$3.996698$   &$3.188049$ &$ 3.096013$ & $ 3.086052$ & $3.083074$& $3.085691$&$ 3.093229$ &$3.137773$ &$3.262639$  & $3.972075$   \\
\hline
$\Exp\tau \,\,\mbox{EXACT}$& $5.639483$& $4.019943$& $3.193551$&$3.10129$ &$3.09131$ & $3.088277$&$3.090955$ &$3.098411$& $3.143053$& $3.269103$  & $4.118051$ \\
\hline
\end{tabular}%
}

\vspace{3mm}
\caption{Expectation $\Exp\tau$ for different $\lambda$ and $u=1$, $x_R=0$.}
\label{Tab:1}
\end{table}

\noindent
Notice that the estimated values of $\Exp \tau$ computed by formula (\ref{pt}) confirm
that the minimum value of the expectation is attained at approximately 
$\lambda^*=1.269812$ (minimum derived from the exact formula \ref{exfor}). Moreover, our computations validate the usefulness of
formulas (\ref{pt}) and (\ref{mr2}).

\begin{figure}[H]
\centering
\begin{minipage}{0.3\textwidth}
\includegraphics[width=\linewidth]{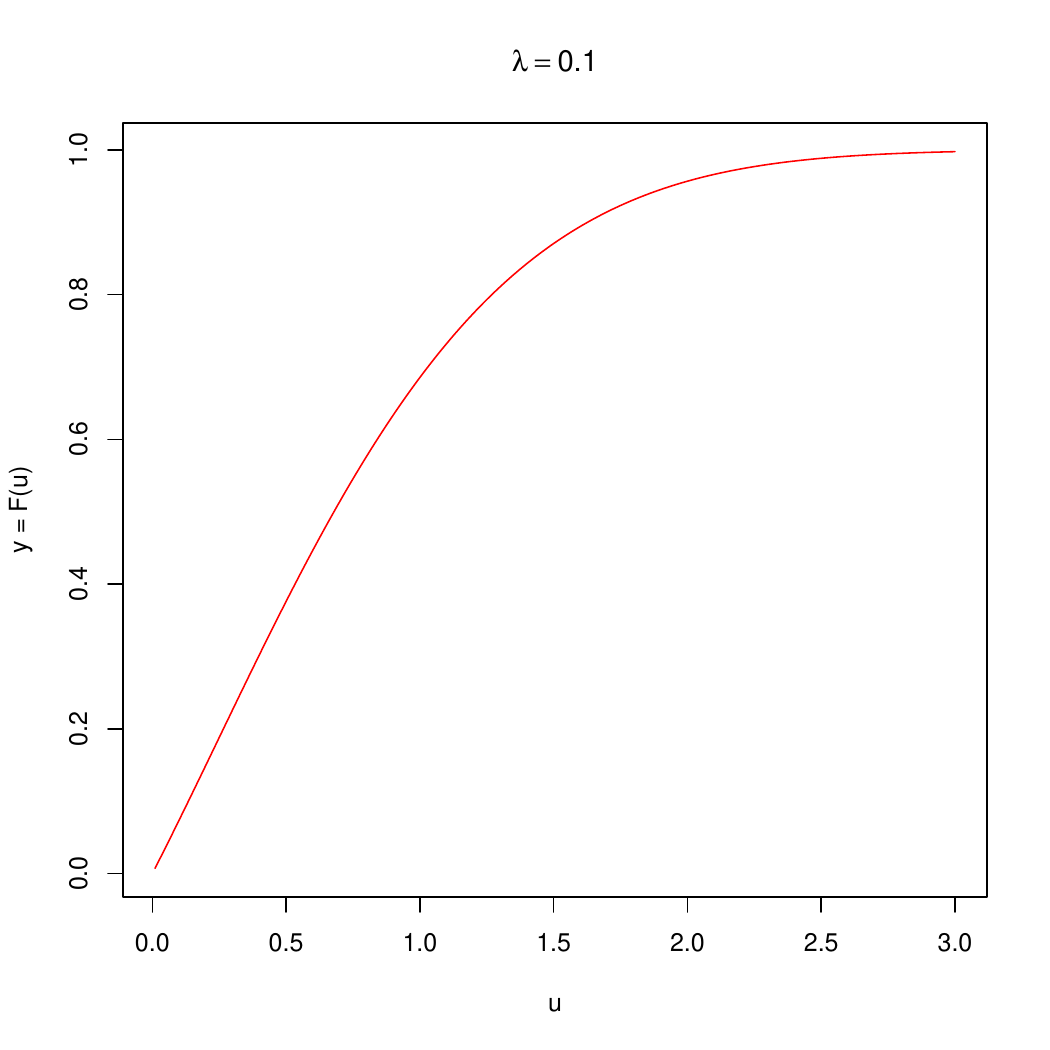}
\end{minipage}%
\begin{minipage}{0.3\textwidth}
\includegraphics[width=\linewidth]{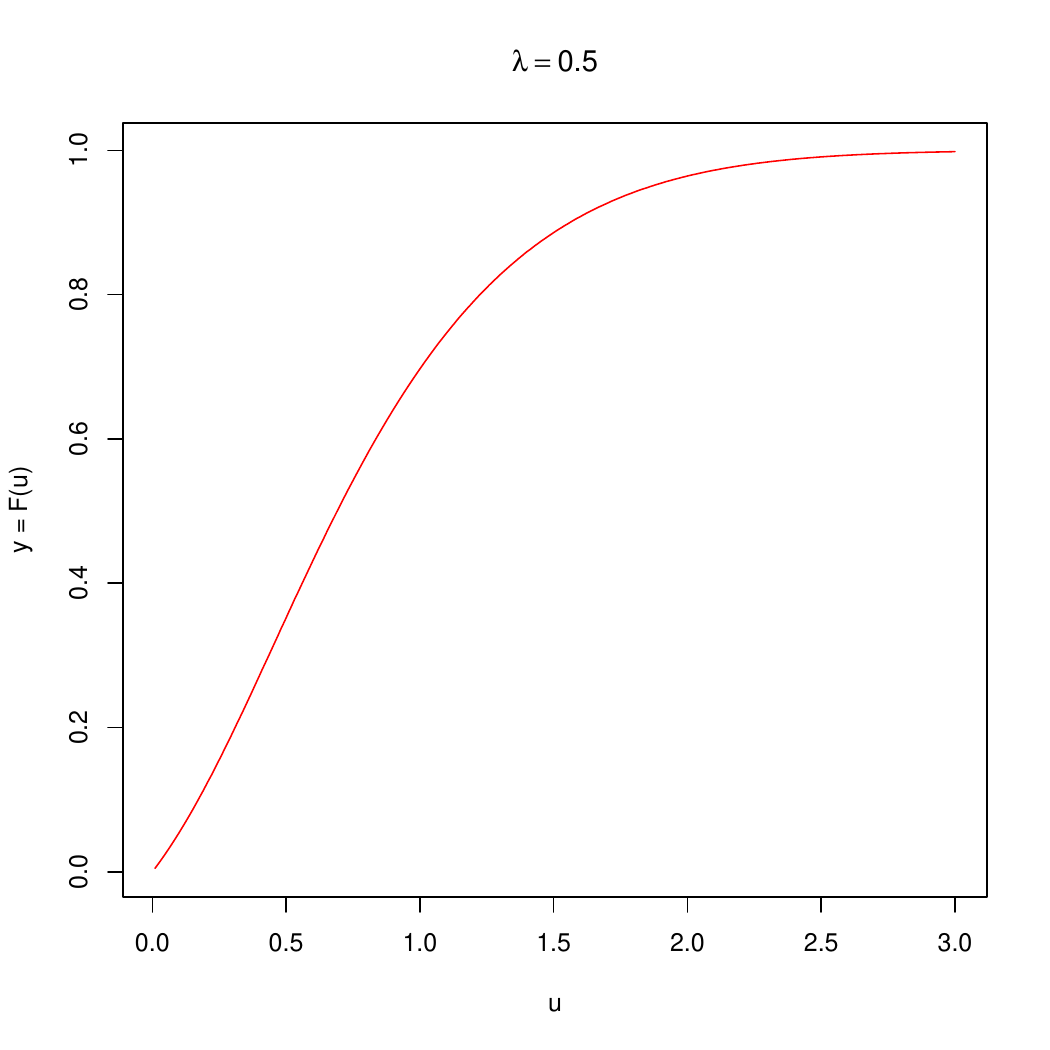}
\end{minipage}%
\begin{minipage}{0.3\textwidth}
\includegraphics[width=\linewidth]{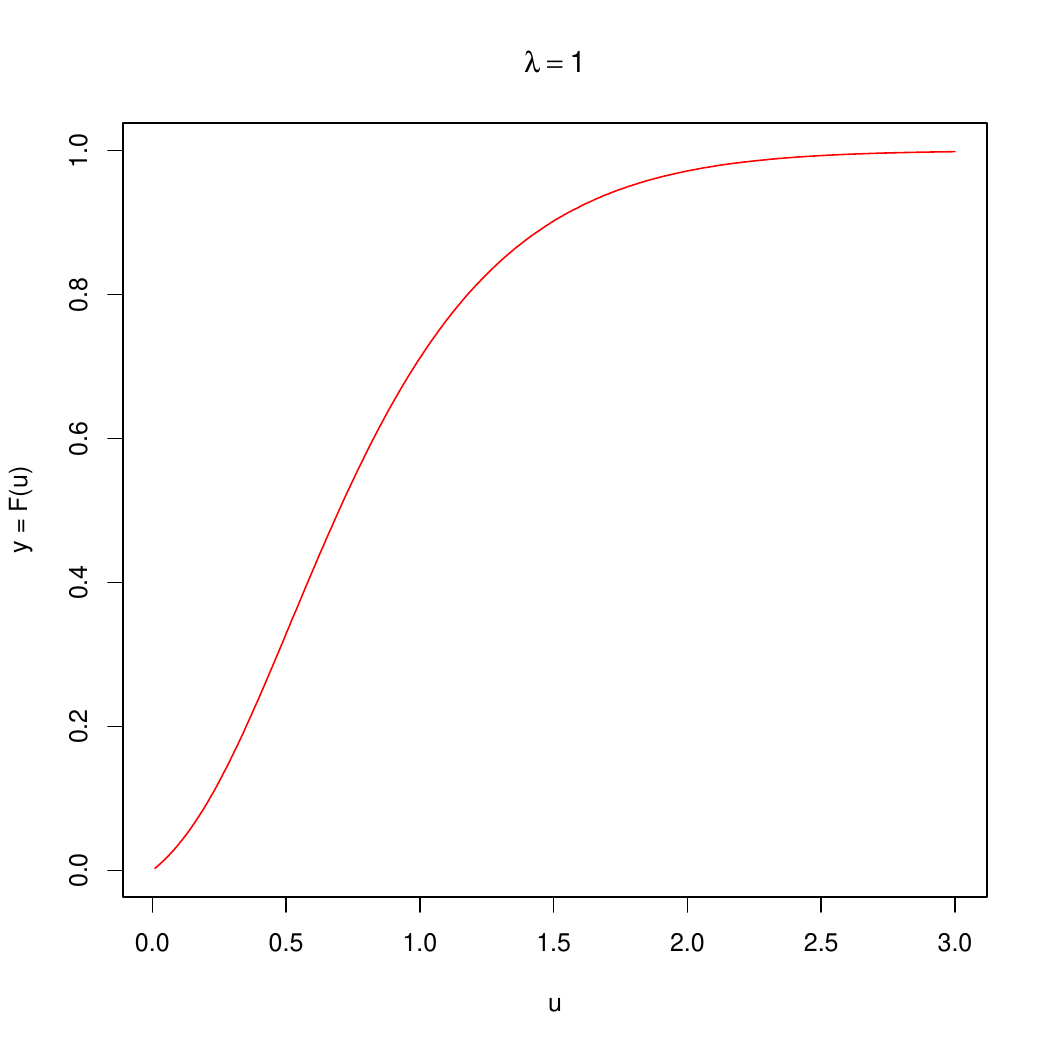}
\end{minipage}

\bigskip

\begin{minipage}{0.3\textwidth}
\includegraphics[width=\linewidth]{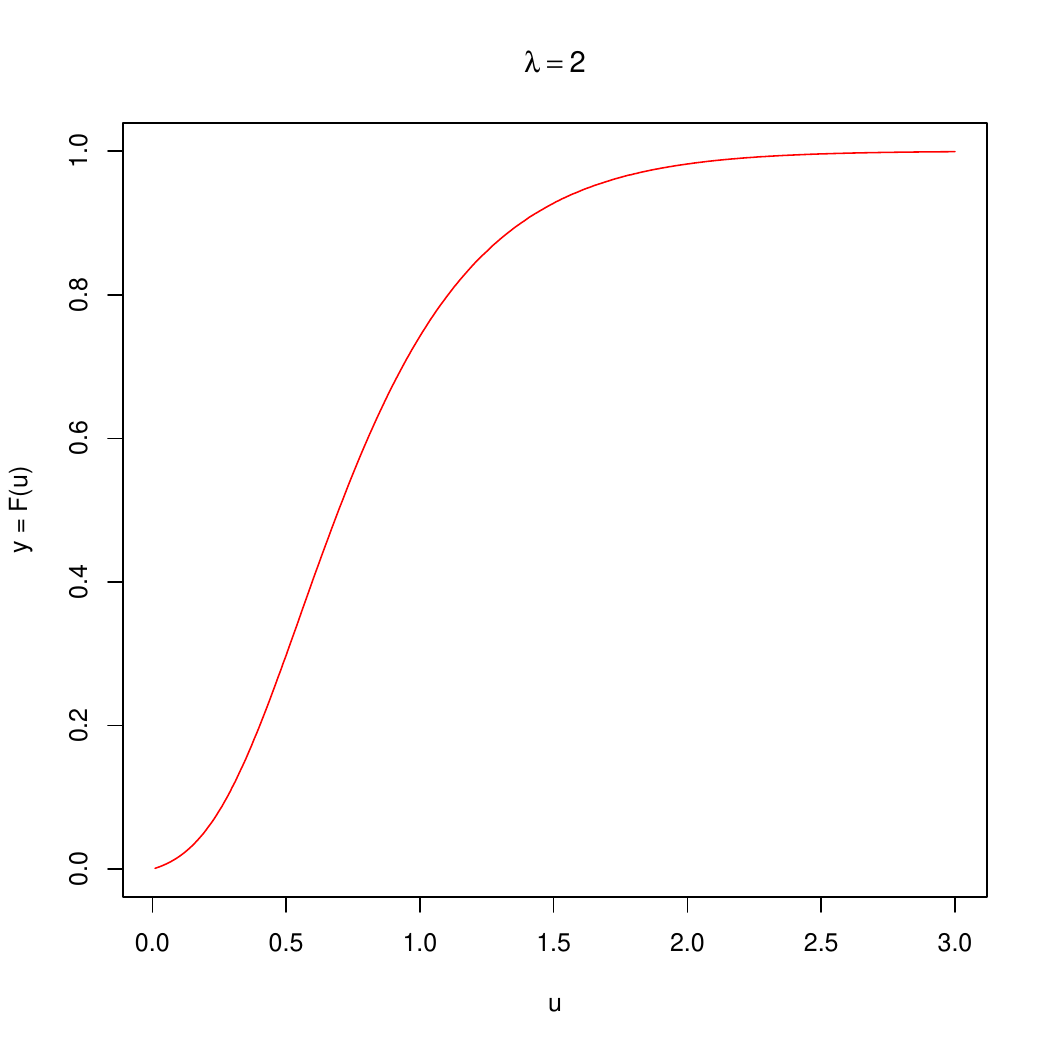}
\end{minipage}%
\begin{minipage}{0.3\textwidth}
\includegraphics[width=\linewidth]{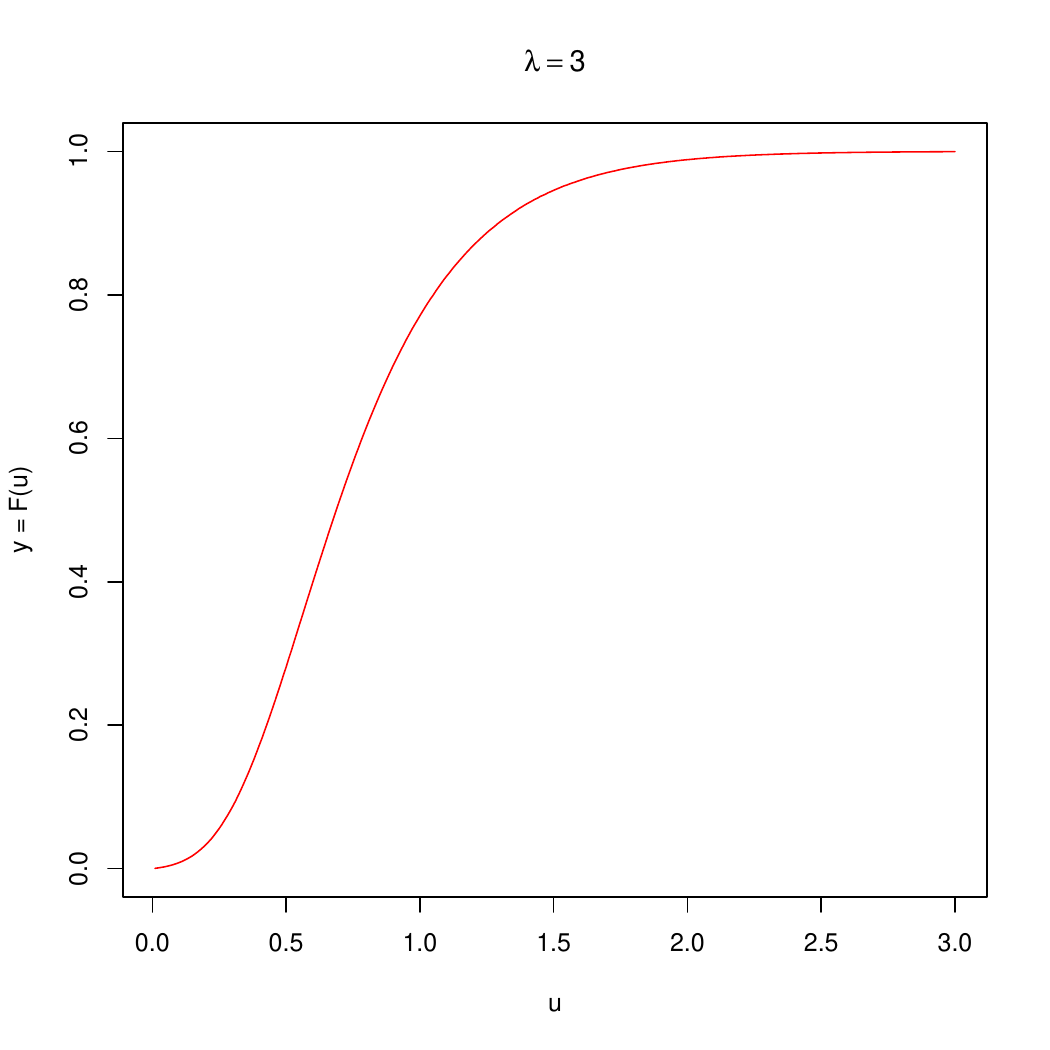}
\end{minipage}%
\begin{minipage}{0.3\textwidth}
\includegraphics[width=\linewidth]{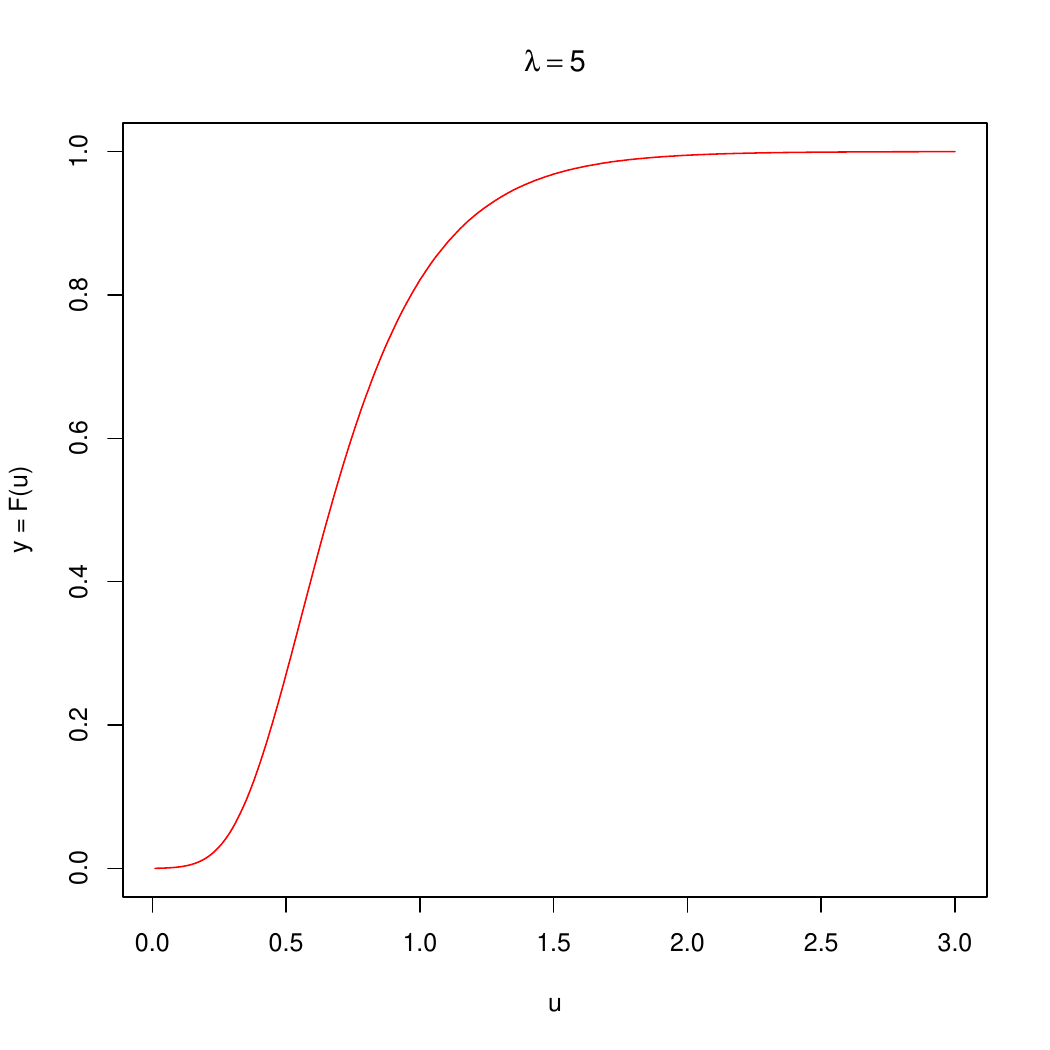}
\end{minipage}

\caption{Cumulative distribution function of $\sup_{t\in [0,1]} X_t$ for $x_R=0$ and $\lambda=0.1, 0.5, 1, 2, 3, 5$.}\label{fig:1}
\end{figure}

Next, using (\ref{mr2}) we approximate values of the cumulative distribution function of $\sup_{t\in[0,1]} X_t$ for $x_R=0$ and $\lambda=0.1, 0.5, 1, 2, 3, 5$, see Fig. \ref{fig:1}. We take $n=100$ summands from (\ref{mr2}) and compute the integrals on
the simplex using the Monte Carlo method with $M=6000$. The discretisation grid of the variable $x$ of the cumulative distribution function $y=F(x)$ is $e=0.01$.

Finally, we numerically approximate the tail distribution of  supremum of a stationary Brownian motion with resetting using
Monte Carlo simulation method and compare it with its asymptotic given in (\ref{asimpBMS}).
For different levels $u$ we compute $\Prob(\sup_{t\in [0,1]}Y_t>u)$ using Monte Carlo simulation method
generating $N=20000$ trajectories with discretization step $e=0.0001$ (MCM in Tab. \ref{Tab:2} and \ref{Tab:3}) and 
compare with the asymptotic given in (\ref{asimpBMS}) (ASYM in Tab. \ref{Tab:2} and \ref{Tab:3}). The probability values for the Monte Carlo method are given with the half-width of the 95\% confidence intervals. The value RATIO in Tab. \ref{Tab:2} and \ref{Tab:3} is the ratio of the Monte Carlo estimation probability
and the asymptotic probability value.
\begin{table}[H]
\centering
\resizebox{\textwidth}{!}{%
\begin{tabular}{l|l|l|l|l|l}
\hline
$u$& $2.5$ & $3$& $3.5$& $4$& $4.5$         \\
\hline
$\Prob(\sup_{t\in [0,1]}Y_t>u)\,\, \mbox{MCM}$&$0.13215\pm 0.00469$ & $0.0501\pm 0.00302$& $0.0186\pm  0.00187 $& $0.0069\pm 0.00115 $&$0.00265\pm 0.000712$ \\
\hline
$\Prob(\sup_{t\in [0,1]}Y_t>u)\,\, \mbox{ASYM}$ & $0.13677$ & $0.05031518$&$0.01850992$   &$0.006809419$ &$0.002505$   \\
\hline
\mbox{RATIO}&$0.96621$ &$0.99572$& $1.0049$& $1.0133$&$ 1.0579$ \\
\hline
\end{tabular}%
}

\vspace{3mm}
\caption{$\Prob(\sup_{t\in [0,1]}Y_t>u)$ for different values of level $u$, $\lambda=2$, and $x_R=1$.}
\label{Tab:2}
\end{table}

\begin{table}[H]
\centering
\resizebox{\textwidth}{!}{%
\begin{tabular}{l|l|l|l|l|l}
\hline
$u$& $2$ & $2.5$& $3$& $3.5$& $4$         \\
\hline
$\Prob(\sup_{t\in [0,1]}Y_t>u)\,\, \mbox{MCM}$&$0.29895\pm 0.00634 $ & $0.0954\pm 0.00407$& $0.0293\pm 0.00234$& $0.0087\pm 0.0013$&$0.00265\pm 0.000712$ \\
\hline
$\Prob(\sup_{t\in [0,1]}Y_t>u)\,\, \mbox{ASYM}$ & $0.3237$ & $0.095115$&$0.027948$   &$0.008212$ &$0.002413$   \\
\hline
\mbox{RATIO}&$0.92353$ &$1.003$& $1.050788$& $1.08673$&$ 1.0982$ \\
\hline
\end{tabular}%
}

\vspace{3mm}
\caption{$\Prob(\sup_{t\in [0,1]}Y_t>u)$ for different values of level $u$, $\lambda=3$, and $x_R=1$.}
\label{Tab:3}
\end{table}

\section{Appendix}\label{app}
In this section, we give complete proofs of all the results presented in Sections \ref{sec2}-\ref{sec4}.
\subsection{Proof of Theorem \ref{jointdist_drift}}\label{app.th1}
The idea of the proof is based on conditioning on $(N_s=k, N_t=l)$ and
considering the last resetting before the epoch $s$ and the last resetting
before the epoch $t$. Without loss of generality, we assume below that $x_0=0.$

First, let us note that for $0\leq s<t$ and $k,l\geq 1$
\[
\Prob(N_s=k, N_t=l)=
\begin{cases}
 e^{-\lambda s}\dfrac{(\lambda s)^k}{k!}\,
 e^{-\lambda (t-s)}\dfrac{[\lambda (t-s)]^{l-k}}{(l-k)!},
 & 0\leq k\leq l,\\[2mm]
 0, & \text{otherwise}.
\end{cases}
\]

Thus, conditioning under $(N_s=k, N_t=l)$, $0\leq k\leq l$, we obtain
\begin{align*}
    \Prob(X^{(c)}_s\leq u, X^{(c)}_t\leq w)
    =&\sum_{0\leq k \leq l}
    \Prob(X^{(c)}_s\leq u, X^{(c)}_t\leq w\mid N_s=k, N_t=l)\Prob(N_s=k, N_t=l)\\
    =&\,e^{-\lambda t}\Prob(X^{(c)}_s\leq u, X^{(c)}_t\leq w\mid N_s=0, N_t=0)\qquad (A)\\
    &\,\,+\,\,\sum_{l=1}^\infty
    e^{-\lambda t}\frac{[\lambda (t-s)]^{l}}{l!}
    \Prob(X^{(c)}_s\leq u, X^{(c)}_t\leq w\mid N_s=0, N_t=l)\qquad (B)\\
    &\,\,+\,\,\sum_{k=1}^\infty
    e^{-\lambda t}\frac{(\lambda s)^{k}}{k!}
    \Prob(X^{(c)}_s\leq u, X^{(c)}_t\leq w\mid N_s=N_t=k)\qquad (C)\\
    &\,\,+\,\,\sum_{\substack{k<l \\ k,l\geq 1}}^\infty
    e^{-\lambda t}\frac{(\lambda s)^{k}}{k!}
    \frac{[\lambda (t-s)]^{l-k}}{(l-k)!}
    \Prob(X^{(c)}_s\leq u, X^{(c)}_t\leq w\mid N_s=k, N_t=l). \qquad (D)
\end{align*}

Recall that the epochs of the jumps of a Poisson process in $[0,t]$ under
condition $N(t)=n$, denoted by $U^{(1)},U^{(2)},\ldots,U^{(n)}$, are the
order statistics of $n$ independent uniformly distributed random variables on $(0,t)$.
The density of $U^{(n)}$ equals $f(x)=\frac{n x^{n-1}}{t^n}$ on $(0,t)$.

Hence,
$(A)=e^{-\lambda t}
\Prob(W^{(c)}_s\leq u, W^{(c)}_t\leq w)$.
The second term $(B)$ equals
\begin{align*}
(B)
&= \Prob(W^{(c)}_s\leq u)
\sum_{l=1}^\infty
e^{-\lambda t}\frac{[\lambda (t-s)]^{l}}{l!}
\Prob(W^{(c)}_{t-s-U^{(l)}}\leq w-x_R)\\
&= \Prob(W^{(c)}_s\leq u)
\sum_{l=1}^\infty
e^{-\lambda t}\frac{[\lambda (t-s)]^{l}}{l!}
\int_0^{t-s}
\Prob(W^{(c)}_{t-s-x}\leq w-x_R)
\frac{l x^{l-1}}{(t-s)^l}\di x\\
&= e^{-\lambda s}\Prob(W^{(c)}_s\leq u)
\lambda\int_0^{t-s}
\Prob(W^{(c)}_x\leq w-x_R)e^{-\lambda x}\di x .
\end{align*}

Let us compute the term $(C)$:
\begin{align*}
(C)
&=\sum_{k=1}^\infty
e^{-\lambda t}\frac{(\lambda s)^k}{k!}
\Prob(W^{(c)}_{s-U^{(k)}}\leq u-x_R,
      W^{(c)}_{t-U^{(k)}}\leq w-x_R)\\
&=e^{-\lambda (t-s)}\lambda
\int_0^s
\Prob(W^{(c)}_{x}\leq u-x_R,
      W^{(c)}_{t-s+x}\leq w-x_R)
e^{-\lambda x}\di x .
\end{align*}
Finally, let us calculate the fourth term $(D)$. Notice that the epochs of jumps on $(0,s)$ are independent of the epochs of jumps on $(s, t]$ under the condition  $(N_s=k,\,N_t=l)$. Thus 
\begin{align*}
    (D)=&\,\, e^{-\lambda t}\sum_{\substack{k<l \\ k,l\geq 1}}^\infty\frac{(\lambda s)^{k}}{k!}\frac{[\lambda (t-s)]^{l-k}}{(l-k)!}\Prob(X^{(c)}_s\leq u, X^{(c)}_t\leq w| N_s=k, N_t-N_s=l-k)\\
    &\,\,=e^{-\lambda t}\sum_{\substack{k<l \\ k,l\geq 1}}^\infty\frac{(\lambda s)^{k}}{k!}\frac{[\lambda (t-s)]^{l-k}}{(l-k)!}\Prob(W^{(c)}_{s-U^{(k)}}\leq u-x_R)
    \Prob(W^{(c)}_{t-s-U^{(l-k)}}\leq w-x_R)\\
    &\,\,=e^{-\lambda t}\sum_{\substack{k<l \\ k,l\geq 1}}^\infty\frac{(\lambda s)^{k}}{k!}\frac{[\lambda (t-s)]^{l-k}}{(l-k)!}
    \int_0^s\Prob(W^{(c)}_{s-x}\leq u-x_R)\frac{k x^{k-1}}{s^k}\di x\\
    &\qquad \times\int_0^{t-s}\Prob(W^{(c)}_{t-s-y}\leq w-x_R)\frac{(l-k)y^{l-k-1}}{(t-s)^{l-k}}\di y\\
    &\,\,=e^{-\lambda t}\sum_{k=1}^\infty\frac{\lambda^{k}}{(k-1)!}
    \int_0^s\Prob(W^{(c)}_x\leq u-x_R)(s-x)^{k-1}\di x\\
    &\qquad\times\sum_{l=k+1}^\infty\frac{\lambda^{l-k}}{(l-k-1)!}
    \int_0^{t-s}\Prob(W^{(c)}_y\leq w-x_R)(t-s-y)^{l-k-1}\di y\\
    &\,\,=e^{-\lambda t}\sum_{k=1}^\infty\frac{\lambda^{k}}{(k-1)!}
    \int_0^s\Prob(W^{(c)}_x\leq u-x_R)(s-x)^{k-1}\di x\\
    &\qquad\times\lambda\int_0^{t-s}\Prob(W^{(c)}_y\leq w-x_R)e^{\lambda(t-s-y)}\di y\\
    &\,\,=\lambda\int_0^s\Prob(W^{(c)}_x\leq u-x_R)e^{-\lambda x}\di x
    \cdot\lambda\int_0^{t-s}\Prob(W^{(c)}_y\leq w-x_R)e^{-\lambda y}\di y\,.
\end{align*}

Collecting $(A)$-$(D)$ completes the proof.
\qed
\subsection{Proof of Theorem \ref{th2}}\label{a.2}
Let us consider the process $X^{(c)}_t$ under the condition that $N_T=n$, $n=0,1,2\ldots$. First, recall that
the epochs of the jumps of a Poisson process in the time interval $[0,T]$ under the condition $N_T=n$ are $n$ order statistics of uniform distribution in $(0,T)$. If $U_1, U_2,\ldots, U_n$ are uniformly distributed in $(0,T)$ then the order statistics form a vector $(U^{(1)}, U^{(2)}, \ldots, U^{(n)})$ with the uniform distribution on
$0<u_1<u_2<\ldots<u_n<T$ (a density function equals $n!/T^n$ in the standard ordered simplex).
Moreover, the distances between adjacent order statistics $D_i=U^{(i)}-U^{(i-1)}$, $i=1,2,\ldots,n$, $U^{(0)}=0$ form a vector $D=(D_1, D_2,\ldots, D_n)$ with a uniform distribution on the standard simplex $u_1+u_2+\ldots+u_n<T$, $u_i>0$, $i=1,2,\ldots,n$.
Note that by the independence of the Wiener process $W$ and the Poisson process $N$ the vector $D$ is independent of the Wiener process $W$.
Let us find the supremum distribution of $X^{(c)}_t$ under the condition that $N_T=n$, $n\geq 1$. By the independence of increments of the Wiener process, we get the following (see Fig. \ref{Pic:1}) 
\begin{align}
\lefteqn{\sup_{t\in [0, T]}  X^{(c)}_t\overset{d}{=}}\label{underNn}\\
&\max\{\sup_{t \in [0,D_1]}W_t^{(c),1}\zm{+x_0}, \sup_{t\in [0,D_2]}W_t^{(c),2}+x_R,\dots, \sup_{t \in [0,D_n]} W_t^{(c),n}+x_R,
\!\sup_{t \in [0, T-D_1-\ldots-D_n]}\!\!\!W_t^{(c),n+1}+x_R\},\nonumber
\end{align}
where $W^{(c),i}$, $i=1,2,\ldots,n+1$  are independent copies of   $W^{(c)}$ being further independent of the vector $D$. Under the condition that $N_T=0$ we have that
$$
\sup_{t\in [0, T]} X^{(c)}_t = \sup_{t\in [0, T]} W^{(c)}_t\zm{+x_0}\,.
$$
Thus, by (\ref{underNn}) we are able to compute the supremum distribution of $X$ conditioning on $N_T=n$, $n\geq 1$ for $u\in \R$
\begin{align*}
\lefteqn{\Prob(\sup_{t\in [0, T]} X^{(c)}_t\leq u| N_T=n)=}\\
&\Prob(\max\{\sup_{t \in [0,D_1]}W_t^{(c),1}\zm{+x_0}, \sup_{t\in [0,D_2]}W_t^{(c),2}+x_R,\dots, \sup_{t \in [0,D_n]} W_t^{(c),n}+x_R,
\sup_{t \in [0, T-D_1-\ldots-D_n]}W_t^{(c),n+1}+x_R\}\leq u)\\
&=\,\Prob(\sup_{t \in [0,D_1]}W_t^{(c),1}\zm{+x_0}\leq u, \sup_{t\in [0,D_2]}W_t^{(c),2}+x_R\leq u,\dots, \\
&\,\,\,\,\,\,\,\,\,\,\,\,\,\,\,\,\,\,\,\,\,\,\,\,\,\,\,\,\,\,\,
\sup_{t\in [0,D_n]} W_t^{(c),n}+x_R\leq u,
\sup_{t \in [0,T-D_1-\ldots-D_n]}W_t^{(c),n+1}+x_R\leq u)\\
&=\,\Exp_D\left[\Prob(\sup_{t \in [0, D_1]}W_t^{(c),1}\zm{+x_0}\leq u|D)\Prob(\sup_{t \in [0,D_2]}W_t^{(c),2}+x_R\leq u|D)\dots\Prob(\sup_{t \in [0,D_n]} W_t^{(c),n}+x_R\leq u|D)\right.\\
&\,\,\,\,\,\,\,\,\,\,\,\,\,\,\,\,\,\,\,\,\,\,\,\,\,\,\,\,\,\,\,\times\left.\Prob(\sup_{t\in [0, T-D_1-\ldots-D_n]}W_t^{(c),n+1}+x_R\leq u|D)\right]\\
&=\,\frac{n!}{T^n} \int\limits_{\substack{s_1 + s_2 + \ldots + s_n < T \\ s_i > 0}} 
\!\!\!\!F(u\zm{-x_0},s_1) 
\prod_{k=2}^n F(u-x_R, s_k)  F(u-x_R, T-s_1-\ldots-s_n)
\, \di s_1 \di s_2 \ldots \di s_n\,.
\end{align*}
Using the total probability formula, we arrive at \eqref{mr2}.
\qed
\subsection{Proof of Theorem \ref{th.sup.as} }\label{a.th.sup.as}
{\underline{Case $x_R\le0$.}} Take $u>0$ and $u>x_R$ and put
\[
p_t(u)\coloneqq 1-F(u,t)
=
\Psi(z_t)
+
e^{-2cu}\Psi(\tilde z_t),
\qquad 
z_t=\frac{u+ct}{\sqrt t}, 
\quad 
\tilde z_t=\frac{u-ct}{\sqrt t},
\qquad t>0.
\]
Fix $\varepsilon\in(0,T]$. For $t\in[\varepsilon,T]$ we have
\[
\inf_{t\in[\varepsilon,T]}\min(z_t,\tilde z_t)
\ge
\frac{u-|c|T}{\sqrt T}
\;\xrightarrow[u\to\infty]{}\;\infty.
\]
Hence, by Mills' ratio, uniformly on $t\in[\varepsilon,T]$,
\[
\Psi(z_t)=\frac{\varphi(z_t)}{z_t}\bigl(1+O(z_t^{-2})\bigr),
\qquad
\Psi(\tilde z_t)=\frac{\varphi(\tilde z_t)}{\tilde z_t}\bigl(1+O(\tilde z_t^{-2})\bigr).
\]
Since $z_t\asymp u$ and $\tilde z_t\asymp u$ uniformly {for  $t\in [\varepsilon,T]$}
{($f(u)\asymp g(u)$ means $C_1 g(u)\leq f(u)\leq C_2 g(u)$ for $C_1>0$, $C_2>0$ and sufficienlty large $u$)},
the remainders may be written as $O(u^{-2})$ uniformly on this interval.
Moreover, the identity
\[
e^{-2cu}\varphi(\tilde z_t)=\varphi(z_t)
\]
holds for all $u$ and $t$, and
\[
\frac{z_t}{\tilde z_t}
=
\frac{u+ct}{u-ct}
=
1+O(u^{-1})
\]
uniformly on $t\in[\varepsilon,T]$.
Consequently,
\[
e^{-2cu}\Psi(\tilde z_t)\sim \Psi(z_t)
\]
uniformly on $t\in[\varepsilon,T]$ as $u\rightarrow\infty$, and plugging this into the exact formula yields
\begin{equation}
p_t(u)\sim 2\Psi\!\left(\frac{u+ct}{\sqrt t}\right),
\qquad u\to\infty,
\label{uniT}
\end{equation}
uniformly on $t\in[\varepsilon,T]$.

By (\ref{ineqF}), we have
    $$
e^{-\lambda T}p_T(u)\leq \Prob(\sup_{t\in [0, T]} X^{(c)}_t > u)\leq 1-(1-p_T(u))e^{-\lambda Tp_T(u-x_R)}\,.
    $$
Using that $e^x=1+x+o(x)$ as $x\rightarrow 0$, we get
$$
e^{-\lambda T}p_T(u)\leq \Prob(\sup_{t\in [0, T]} X^{(c)}_t > u)\leq \lambda Tp_T(u-x_R)+p_T(u)e^{-\lambda Tp_T(u-x_R)}+o(p_T(u-x_R)).
$$
Hence from \eqref{uniT},  uniformly on $t\in (\epsilon,T]$
\begin{eqnarray*}
\lim_{u\rightarrow\infty}\frac{p_t(u-x_R)}{p_t(u)}=\left\{
\begin{array}{ll}
 0    & \mbox{if} \,\,\,x_R<0 \\
 1   & \mbox{if} \,\,\,x_R=0 
\end{array}
\right.
\end{eqnarray*} 
leads to
\begin{eqnarray}
2e^{-\lambda T}\leq \liminf_{u\rightarrow\infty}\frac{\Prob(\sup_{t\in [0, T]} X^{(c)}_t > u)}{\Psi\left(\frac{u+cT}{\sqrt{T}}\right)}
\leq\limsup_{u\rightarrow\infty}\frac{\Prob(\sup_{t\in [0, T]} X^{(c)}_t > u)}{\Psi\left(\frac{u+cT}{\sqrt{T}}\right)}\leq 2(\lambda T+1)\, .   
\label{inq.Z}
\end{eqnarray}
In view of (\ref{inq.Z}),
in order to derive the exact asymptotics, as claimed in $i)$, it suffices to focus on  the upper bound
for $\Prob(\sup_{t\in [0, T]} X^{(c)}_t > u)$.
Let $\tau$ be the first reset time of $X^{(c)}_t, t\ge0$, and 
$\tilde{X}^{(c)}_t,t\ge0$ be the resetting process with the resetting
value $\tilde{x}_R=0$, which is independent of ${X}^{(c)}_t,t\ge0$. 
For $\varepsilon\in(0,T/2)$, we have
\begin{eqnarray}
 \Prob(\sup_{t\in [0, T]} X^{(c)}_t > u)&\le& 
    \Prob(\sup_{t\in [0, T]} X^{(c)}_t > u,\tau\in[0,\varepsilon]\cup [T-\varepsilon,T])
    +\Prob(\sup_{t\in [0, T]} X^{(c)}_t > u,\tau\in[\varepsilon,T-\varepsilon])\nonumber\\
   &&
   + \Prob(\sup_{t\in [0, T]} X^{(c)}_t > u,\tau>T)
   =:
   P_1(u)+P_2(u)+P_3(u).
\end{eqnarray}
We analyze each component of the above display separately, beginning with $P_3(u)$, for which we have
\begin{eqnarray*}
P_3(u)=\Prob(\tau>T)\Prob(\sup_{t\in [0, T]} (W_t-ct) > u)\sim 2e^{-\lambda T}\Psi\left(\frac{u+cT}{\sqrt{T}}\right),
\end{eqnarray*}
as $u\to\infty$. 
Next,  we have
\begin{eqnarray}
P_1(u)&\le&
\Prob(\sup_{t\in [0, \tau]} X^{(c)}_t > u,\tau\in[0,\varepsilon]\cup [T-\varepsilon,T])
+
\Prob(\sup_{t\in [\tau, T]} X^{(c)}_t > u,\tau\in[0,\varepsilon]\cup [T-\varepsilon,T]). \nonumber
\end{eqnarray}
Using that
\begin{eqnarray}
    \Prob(\sup_{t\in [0, \tau]} X^{(c)}_t > u,\tau\in[0,\varepsilon]\cup [T-\varepsilon,T])
&=&
\Prob(\sup_{t\in [0, \tau]} (W_t-ct) > u,\tau\in[0,\varepsilon]\cup [T-\varepsilon,T])\nonumber\\
&\le&
\Prob(\sup_{t\in [0, T]} (W_t-ct) > u)\Prob(\tau\in[0,\varepsilon]\cup [T-\varepsilon,T])\nonumber\\
&\le&
2\varepsilon \lambda\Prob(\sup_{t\in [0, T]} (W_t-ct) > u)\nonumber
\end{eqnarray}
and
\begin{eqnarray}
\Prob(\sup_{t\in [\tau, T]} X^{(c)}_t > u,\tau\in[0,\varepsilon]\cup [T-\varepsilon,T])
&=&
\Prob(\sup_{t\in [0, T-\tau]} \tilde{X}^{(c)}_t+x_R > u,\tau\in[0,\varepsilon]\cup [T-\varepsilon,T])\nonumber\\
&\le&
\Prob(\sup_{t\in [0, T]} \tilde{X}^{(c)}_t > u)\Prob(\tau\in[0,\varepsilon]\cup [T-\varepsilon,T])
\nonumber\\
&\le&
2\varepsilon \lambda\Prob(\sup_{t\in [0, T]} \tilde{X}^{(c)}_t > u),\nonumber
\end{eqnarray}
in view of (\ref{uniT}) and (\ref{inq.Z}),
we deduce that, as $u\to\infty$,
\begin{eqnarray}
    \Prob(\sup_{t\in [0, \tau]} X^{(c)}_t > u,\tau\in[0,\varepsilon]\cup [T-\varepsilon,T])
&\le& 
2\varepsilon \lambda(3+2\lambda T) \Psi\left(\frac{u+cT}{\sqrt{T}}\right) .\label{as.in}
\end{eqnarray}
Finally,
\begin{eqnarray*}
P_2(u)&\le&
\Prob(\sup_{t\in [0, \tau]} X^{(c)}_t > u,\tau\in[\varepsilon,T-\varepsilon])
+
\Prob(\sup_{t\in [\tau, T]} X^{(c)}_t > u,\tau\in[\varepsilon,T-\varepsilon])\\
&\le&
\Prob(\sup_{t\in [0, T-\varepsilon]} (W_t-ct) > u)
+
\Prob(\sup_{t\in [0, T-\tau]} \tilde{X}^{(c)}_t+x_R > u,\tau\in[\varepsilon,T-\varepsilon])\nonumber\\
&\le&
\Prob(\sup_{t\in [0, T-\varepsilon]} (W_t-ct) > u)
+
\Prob(\sup_{t\in [0, T-\varepsilon]} \tilde{X}^{(c)}_t > u)
=
o\left(\Psi\left(\frac{u+cT}{\sqrt{T}}\right)\right),
\end{eqnarray*}
as $u\to\infty$,
by the combination of (\ref{uniT}) with (\ref{inq.Z}).
Thus, 
\[
\limsup_{u\to\infty}\frac{\Prob(\sup_{t\in [0, T]} X^{(c)}_t > u)}{\Psi\left(\frac{u+cT}{\sqrt{T}}\right)}\le
2e^{-\lambda T},
\]
which combined with
(\ref{inq.Z}) completes the proof of case $i)$.
\\
{\underline{Case $x_R>0$.}}
With the same notation as in case $i)$, we begin with the observation that
\begin{eqnarray} 
    \Prob(\sup_{t\in [\tau, T]} X^{(c)}_t > u, \tau \le \log^2 (u) u^{-2} )
    &\le&
    \Prob(\sup_{t\in [0, T]} X^{(c)}_t > u)
    \le
    \Prob(\sup_{t\in [0, T]} X^{(c)}_t > u, \tau \le \log^2(u) u^{-2} )\nonumber\\
   &+&
    \Prob(\sup_{t\in [0, T]} X^{(c)}_t > u, \tau > \log^2(u)  u^{-2} ) =: R_1(u)+R_2(u)\,.\qquad\label{bounds.2}
\end{eqnarray}
We note that $R_1(u)\zm{\le}R_{11}(u)+R_{12}(u)$, where
\begin{eqnarray}
R_{11}(u)\coloneqq\Prob(\sup_{t\in [0, \tau]} X^{(c)}_t > u, \tau \le \log^2  (u) u^{-2} )
 &=&
    \Prob(\sup_{t\in [0, \tau]} (W_t-ct) > u, \tau \le \log^2  (u) u^{-2} )\nonumber\\
  &\le&
     \Prob(\sup_{t\in [0, T]} (W_t-ct) > u)\nonumber\\
  &\sim& 2\Psi\left(\frac{u+cT}{\sqrt{T}}\right)\label{R_11}
 \end{eqnarray}
and 
\begin{eqnarray*}
R_{12}(u):=    \Prob(\sup_{t\in [\tau, T]} X^{(c)}_t > u, \tau \le \log^2  (u) u^{-2} )
    &=&
    \int_{0}^{\log^2  (u) u^{-2}}
         \Prob(\sup_{t\in [0, T-x]} \tilde{X}^{(c)}_t +x_R> u) \lambda e^{-\lambda x} dx.
\end{eqnarray*}
The asymptotics of
\[
\Prob(\sup_{t\in [0, T-x]} \tilde{X}^{(c)}_t> u-x_R)=:q_{T-x}(u).
\]
follows  by item $i)$,   which holds   uniformly  for $x\in [0,\delta], \delta \in [0,T/2)$.
Consequently, we have
\begin{eqnarray}
  R_{12}(u)
&\sim&
2e^{-\lambda T}\int_{0}^{\log^2  (u) u^{-2}}
\Psi\left(\frac{u-x_R+c\zm{(T-x)}}{\sqrt{T-x}}\right)
         \lambda e^{-\lambda x} dx\nonumber\\
&=&         
2e^{-\lambda T}u^{-2}\int_{0}^{\log^2  (u)}
\Psi\left(\frac{u-x_R+c\zm{(T-x/u^2)}}{\sqrt{T-x/u^2}}\right)
         \lambda e^{-\lambda x/u^2} dx\nonumber\\
&\sim&
2 e^{-\lambda T}u^{-2}\int_{0}^{\log^2  (u)}
\Psi\left(\frac{u-x_R+c\zm{(T-x/u^2)}}{\sqrt{T}}\right) e^{-x/(2T^2)}\lambda e^{-\lambda x/u^2}
          dx\nonumber\\
&\sim&
4\lambda T^2 e^{-\lambda T}u^{-2}
\Psi\left(\frac{u-x_R+cT}{\sqrt{T}}\right) ,\label{R_12}
\end{eqnarray}
as $u\to\infty$.
Moreover
\begin{eqnarray}
    R_{2}(u)
    &\le&
    \Prob(\sup_{t\in [0, \tau]} X^{(c)}_t > u, \tau > \log^2  (u) u^{-2} ) 
    +
    \Prob(\sup_{t\in [\tau, T]} X^{(c)}_t > u, \tau > \log^2  (u) u^{-2} ) \nonumber\\
    &\le&
    \Prob(\sup_{t\in [0, T]} (W_t-ct) > u)
    +
    \Prob(\sup_{t\in [0,T-\log^2  (u) u^{-2}]} \tilde{X}^{(c)}_t+x_R > u ) \nonumber\\
    &\sim&
    2  \Psi\left(\frac{u+cT}{\sqrt{T}}\right)
    +
    2e^{-\lambda T}\Psi\left(\frac{u-x_R+cT}{\sqrt{T}}\right) e^{-\frac{\log^2(u)}{2T^2}},\label{R_2}
\end{eqnarray}
as $u\to\infty$.
Combining (\ref{R_11}), (\ref{R_12}), (\ref{R_2}) with (\ref{bounds.2}), we conclude that
\[
 \Prob(\sup_{t\in [0, T]} X^{(c)}_t > u)
   \sim 
   4\lambda T^2 e^{-\lambda T}u^{-2}
\Psi\left(\frac{u-x_R+cT}{\sqrt{T}}\right),
\]
as $u\to\infty$, which completes the proof.
\qed

\subsection{Proof of Theorem \ref{thm:BM-window-reset-general-x0}}\label{a.thm:BM-window-reset-general-x0}
We first prove the corresponding uniform window asymptotic for Brownian motion with drift extending \cite[Thm 2.1]{dkebicki2016parisian}.
For all $0<a<b<\infty$ and $r\ge0$,
\begin{equation}\label{eq:uniform-window-proof}
\lim_{u\to\infty}\sup_{\substack{T\in[a,b]\\ \delta\in[a,b]}}
\left|
\frac{
\Prob\left( \inf_{ s \in [T,T+\delta]} W_s^{(c)} >u,\; W_T^{(c)} > u+\frac{r}{u+cT}\right)
}{
K_{c,\delta}\,L(r/T)\,
\frac{T}{u}\,
\Psi\!\Big(\frac{u+cT}{\sqrt{T}}\Big)
}
-1
\right|
= 0,
\end{equation}
where  
$$  
K_{c,\delta}
=
\frac{2}{\sqrt{\delta}}\,
\varphi\big(c\sqrt{\delta}\big)
-2c \Psi\big( c \sqrt{\delta}\big).
\qquad 
$$
Fix $(T,\delta)\in[a,b]^2$.  
By the independent increments property
\[
\inf_{t\in[T,T+\delta]}(W_t-ct)
=
(W_T-cT)+\inf_{s\in[0,\delta]}(\widetilde W_s-cs),
\]
where $\widetilde W$ is an independent standard Brownian motion. Let
\[
X\coloneqq W_T,\qquad
Y\coloneqq \sup_{s\in[0,\delta]}(\widetilde W_s+cs)
\]
and 
set
\[
v=v(u,T)\coloneqq u+cT+\frac{r}{u+cT}.
\]
With this notation \begin{equation}\label{eq:XYrepr-proof}
G_T(u,\delta)\coloneqq \pk{\inf_{t\in[T,T+\delta]}(W_t-ct)>u,\;W_T>v}
=
\pk{X-Y>u+cT,\ X>v}.
\end{equation}
In view of  (\ref{probF})
\begin{equation}\label{eq:Ycdf-proof}
\pk{Y\le y}
=
\Phi\!\Big(\frac{y-c\delta}{\sqrt{\delta}}\Big)
-
e^{2cy}\,
\Phi\!\Big(\frac{-y-c\delta}{\sqrt{\delta}}\Big),
\qquad y\ge0.
\end{equation}
Hence $Y$ has a continuous density $f_\delta$ on $(0,\infty)$,
continuous at $0$, and differentiating \eqref{eq:Ycdf-proof} gives
\begin{equation}\label{eq:fdelta0-proof}
f_\delta(0)=K_{c,\delta}>0.
\end{equation}
Since $\delta\in[a,b]$ we have further
\begin{equation}\label{eq:fdelta-unif0-proof}
\lim_{y\downarrow0}\ \sup_{\delta\in[a,b]}\,|f_\delta(y)-K_{c,\delta}|=0
\end{equation}
and
\begin{equation}\label{eq:fdelta-bdd-proof}
\sup_{\delta\in[a,b]}\ \sup_{y\ge0} f_\delta(y)<\infty.
\end{equation}
Define
\[
z(u,T)\coloneqq\frac{u+cT}{\sqrt{T}},
\qquad
w(u,T)\coloneqq\frac{u+cT}{T}.
\]   Since $T\in[a,b]$,
\[
z(u,T)\to\infty\quad\text{and}\quad w(u,T)\to\infty
\qquad\text{uniformly as }u\to\infty.
\]
Hence, by the standard Gaussian tail ratio, for every fixed $M>0$
\begin{equation}\label{eq:ratio-unif-proof}
\lim_{u \to \infty} \sup_{T\in[a,b],x\in[0,M] }
\left|
\frac{\pk{X>u+cT+x/w(u,T)}}{\pk{X>u+cT}}-e^{-x}
\right|
= 0.
\end{equation}
Moreover, there exist constants $C,\eta>0$ (depending only on $a,b,c$) such that for all $u$ large enough,
\begin{equation}\label{eq:ratio-dom-proof}
\sup_{T\in[a,b]}
\frac{\pk{X>u+cT+x/w(u,T)}}{\pk{X>u+cT}}
\le
C e^{-\eta x},
\qquad x\ge0.
\end{equation}
From \eqref{eq:XYrepr-proof} and independence of $X$ and $Y$
\[
G_T(u,\delta)
=\int_0^\infty \pk{X>\max(u+cT+y,\ v)}\,f_\delta(y)\,dy=
\frac{\pk{X>u+cT}}{w(u,T)}\,H_{u,T,\delta},
\]
where
\[
H_{u,T,\delta}
\coloneqq
\int_0^\infty
\frac{\pk{X>u+cT+x^*/w(u,T)}}{\pk{X>u+cT}}\,
f_\delta\!\left(\frac{x}{w(u,T)}\right)\,dx,
\qquad
x^*\coloneqq\max\left(x,\frac{r}{T}\right).
\]
We shall show that, uniformly for \((T,\delta)\in[a,b]^2\)
\[
H_{u,T,\delta}
\longrightarrow
K_{c,\delta}\int_0^\infty e^{-x^*}\,dx
=
K_{c,\delta}L(r/T),
\qquad u\to\infty .
\]
On \([0,M], M>0\), by \eqref{eq:ratio-unif-proof} and
\eqref{eq:fdelta-unif0-proof}, the integrand converges uniformly in
\((T,\delta)\in[a,b]^2\) to
$ 
e^{-x^*}K_{c,\delta}$ 
implying 
\[ \lim_{u\to \infty}
\int_0^M
\frac{\pk{X>u+cT+x^*/w(u,T)}}{\pk{X>u+cT}}\,
f_\delta\!\left(\frac{x}{w(u,T)}\right)\,dx
=
K_{c,\delta}\int_0^M e^{-x^*}\,dx
\]
uniformly for \((T,\delta)\in[a,b]^2\).

It remains to control the tail of the normalized integral. By
\eqref{eq:ratio-dom-proof} and \eqref{eq:fdelta-bdd-proof}, there exists
\(C'>0\) such that, for all sufficiently large \(u\)
\[
\sup_{(T,\delta)\in[a,b]^2}
\int_M^\infty
\frac{\pk{X>u+cT+x^*/w(u,T)}}{\pk{X>u+cT}}\,
f_\delta\!\left(\frac{x}{w(u,T)}\right)\,dx
\le
C'\int_M^\infty e^{-\eta x^*}\,dx .
\]
Since \(x^*=\max(x,r/T)\ge x\), the right-hand side is bounded by
\[
C'\int_M^\infty e^{-\eta x}\,dx
=
\frac{C'}{\eta}e^{-\eta M},
\]
which tends to zero as \(M\to\infty\), uniformly in \((T,\delta)\in[a,b]^2\).

Consequently, we obtain
\[\lim_{u \to \infty}
H_{u,T,\delta}
=
K_{c,\delta}\int_0^\infty e^{-x^*}\,dx
=
K_{c,\delta}L(r/T)
\]
uniformly for \((T,\delta)\in[a,b]^2\). Therefore, 
\[
G_T(u,\delta)
\sim
K_{c,\delta}L(r/T)\frac{T}{u}
\Psi\!\left(\frac{u+cT}{\sqrt T}\right),
\]
uniformly for \((T,\delta)\in[a,b]^2\), hence 
\eqref{eq:uniform-window-proof} follows.

\underline{Case $x_0<x_R$.} 
Let $r\ge0$ and
\[
v=u+\frac{r}{u+cT}.
\]
Assume $u>x_R$ and   
recall the  exact decomposition \eqref{eq:inf-window-exact-v} 
\begin{equation}\label{eq:decomp-i}
\begin{aligned}
\pk{\inf_{t\in[T,T+\Delta]}X_t^{(c)} >u,\; X_{T}^{(c)} >v}
&=
e^{-\lambda (T+\Delta)}\,
F_{T+\Delta}^{(v-x_0+cT)}(u-x_0,\Delta)
+e^{-\lambda\Delta} J(u),
\end{aligned}
\end{equation}
where
\[
F_s^{(w)}(a,\Delta)
:=
\pk{\inf_{t\in[s-\Delta,s]}(W_t-ct)>a,\; W_{s-\Delta}>w}, \qquad J(u)= 
 \int_{0}^{T}\lambda e^{-\lambda s}\,
F_{s+\Delta}^{(v-x_R+cs)}(u-x_R,\Delta)ds.
\]
 
Set  $a:=u-x_R$ and 
fix $\varepsilon\in(0,T)$.
For $s\in[T-\varepsilon,T]$, we apply \eqref{eq:uniform-window-proof} with $(T,\delta)=(s,\Delta)$
and level $a=u-x_R$. Note that the constraint level at time $s$ is
\[
W_s > v-x_R+cs = a+cs+\frac{r}{u+cT}
=
a+cs+\frac{r}{a+cs}\,(1+o(1)),
\]
uniformly for $s\in[T-\varepsilon,T]$ (since $a\to\infty$ and $cs$ is bounded).
Thus the same tail-ratio argument used to derive \eqref{eq:uniform-window-proof} yields that, uniformly for
$s\in[T-\varepsilon,T]$
\begin{equation}\label{eq:Fs-asymp}
F_{s+\Delta}^{(v-x_R+cs)}(a,\Delta)
\sim
K_{c,\Delta}\,L(r/s)\,\frac{s}{a}\,
\Psi\!\Big(\frac{a+cs}{\sqrt{s}}\Big),
\qquad u\to\infty.
\end{equation}
Moreover, $L(r/s)=L(r/T)+o(1)$ uniformly on $[T-\varepsilon,T]$.

For $s\in(0,T-\varepsilon]$, the Gaussian rate
\[
I_a(s)\coloneqq \frac{(a+cs)^2}{2s}
\]
satisfies $I_a(s)\ge I_a(T-\varepsilon)$ and $I_a(T-\varepsilon)-I_a(T)\asymp a^2$.
Hence the contribution of $(0,T-\varepsilon]$ to $J(u)$ is exponentially negligible compared to the
contribution from $[T-\varepsilon,T]$ implying as $u\to \infty$  
\begin{equation}\label{eq:J-local}
J(u)
\sim
\int_{T-\varepsilon}^{T}\lambda e^{-\lambda s}\,
K_{c,\Delta}\,L(r/T)\,\frac{s}{a}\,
\Psi\!\Big(\frac{a+cs}{\sqrt{s}}\Big)\,ds.
\end{equation}
 
Set
\[
z_a\coloneqq \frac{a+cT}{\sqrt{T}}\sim \frac{a}{\sqrt{T}},
\qquad
g(s)\coloneqq \frac{a+cs}{\sqrt{s}}.
\]
The minimum of $g(s)$  over $s\in(0,T]$ is attained at $s=T$ for all large $a$,
and the correct boundary-layer scale is $a^{-2}$.
A Taylor expansion yields, locally uniformly for bounded $h$
\begin{equation}\label{eq:g-expand}
g\Big(T-\frac{h}{a^2}\Big)
=
z_a+\frac{h}{2T^2}\,\frac{1}{z_a}+o\!\Big(\frac{1}{z_a}\Big),
\qquad a\to\infty.
\end{equation}
In view of \eqref{eq:ratio-unif-proof}
\[
\lim_{a \to \infty} \frac{\Psi(g(T-h/a^2))}{\Psi(z_a)}
=
\exp\!\Big(-\frac{h}{2T^2}\Big) 
\]
locally uniformly in $h§$. Hence, by dominated convergence applied to \eqref{eq:J-local}, which is justified  from \eqref{eq:uniform-window-proof}, we obtain  as $u \to \infty$
\begin{align}
J(u)
&\sim
\lambda e^{-\lambda T}\,K_{c,\Delta}\,L(r/T)\,\frac{T}{a}\,
\Psi(z_a)\,
\frac{1}{a^2}\int_{0}^{\infty}\exp\!\Big(-\frac{h}{2T^2}\Big)\,dh
\notag\\
&\sim 
2\lambda e^{-\lambda T}\,K_{c,\Delta}\,L(r/T)\,
\frac{T^3}{u^3}\,
\Psi\!\Big(\frac{u-x_R+cT}{\sqrt{T}}\Big).\notag 
\end{align}
Finally, application of the uniform window asymptotic \eqref{eq:uniform-window-proof}  yields
\begin{equation}\label{eq:first-term-i}
F_{T+\Delta}^{(v-x_0+cT)}(u-x_0,\Delta)
\sim
K_{c,\Delta}\,L(r/T)\,\frac{T}{u}\,
\Psi\!\Big(\frac{u-x_0+cT}{\sqrt{T}}\Big)=o\!\big(J(u)\big).
\end{equation}  
Hence the claim follows.

\underline{Case $x_0\ge  x_R$.}   
 From \eqref{eq:first-term-i} and the asymptotics of $J(u)$,   using that $x_0\ge  x_R$, we have   that $J(u)$ is dominated by $F_{T+\Delta}^{(v-x_0+cT)}(u-x_0,\Delta)$ establishing the proof. 
 \qed

\subsection{Proof of Proposition \ref{prop.stat}}\label{s.a.prop.stat.}
First, let us show that the Laplace distribution given in \eqref{finfty} and \eqref{inftydistr}
is an invariant probability measure for the given Markov process. 
Thus, using the fact that the time since the last jump at time $t$ under condition $N_t\geq 1$ is exponentially distributed with parameter $\lambda>0$ on $(0,t)$ we get
\begin{align*}
\Prob(Y^{(c)}_t\leq u)=&\,
e^{-\lambda t}\Prob(Y^{(c)}_t\leq u | N_t=0)+(1-e^{-\lambda t})\Prob(Y^{(c)}_t\leq u | N_t\geq 1)\\
&=e^{-\lambda t} \Prob(Y^{(c)}_0+W^{(c)}_t\leq u)+\lambda\int_0^t\Prob(W^{(c)}_s+x_R\leq u)e^{-\lambda s}\di s\\
&=e^{-\lambda t}\int_0^\infty \di s \lambda e^{-\lambda s}\left[\int_{-\infty}^\infty\Prob(x+W^{(c)}_t\leq u)\di_x\Prob(W^{(c)}_s+x_R\leq x)\right]\\
&\,\,\,\,\,+ \lambda\int_0^t\Prob(W^{(c)}_s+x_R\leq u)e^{-\lambda s}\di s\\
&=e^{-\lambda t}\int_0^\infty \lambda e^{-\lambda s}\Prob(W^{(c)}_{s+t}+x_R\leq u)\di s + \lambda\int_0^t\Prob(W^{(c)}_s+x_R\leq u)e^{-\lambda s}\di s\\
&=\lambda\int_t^\infty \Prob(W^{(c)}_{s}+x_R\leq u)e^{-\lambda s}\di s + \lambda\int_0^t\Prob(W^{(c)}_s+x_R\leq u)e^{-\lambda s}\di s\\
&=\lambda\int_0^\infty\Prob(W^{(c)}_s+x_R\leq u)e^{-\lambda s}\di s =\Prob(Y^{(c)}_0\leq u)\,,
\end{align*}
where in the third and last equalities we used (\ref{inftydistr}), in the third equality we also applied independence of $Y^{(c)}_0$ and $W^{(c)}$ and in the fourth equality we used that fact that the sum of two independent normal random variables with variances $s$ and $t$ is normal with variance $s+t$. So we have $Y^{(c)}_t\stackrel{d}{=}Y^{(c)}_0$ for all $t\geq 0$.
To prove that finite-dimensional distributions of $Y_t$ are invariant under time shifting, let us notice that $Y^{(c)}_t$ is a time-homogeneous Markov process. Let
\begin{equation}\label{kernel}
\Prob_t(x, \di y)=\Prob_x(Y^{(c)}_t\in \di y)
\end{equation}
be the transition kernel where $\Prob_x$ denotes probability under condition $Y^{(c)}_0=x$. As we showed above, the measure $\mu(\di x)=f_{X^{(c)}_\infty}(x)\di x$  
is invariant for this transition kernel that is
$$
\int_{-\infty}^\infty \Prob_t(x, A)\mu(\di x)=\mu(A)
$$
for any Borel measurable set $A\in\RL$.
Let $0\leq t_1<\ldots<t_n$ and $h\geq 0$ and $\Phi: \RL^n\rightarrow\RL$ be a measurable bounded function. Using standard arguments for Markov processes, we easily get
$$
\Exp \Phi(Y^{(c)}_{t_1+h},\ldots, Y^{(c)}_{t_n+h})=
\Exp \Phi(Y^{(c)}_{t_1},\ldots, Y^{(c)}_{t_n})\,,
$$ 
which proves stationarity of the process $Y^{(c)}$.
Clearly, the transition kernel defined in \eqref{kernel} satisfies 
$
\Prob_t(x, \di y)=\Prob_x(X^{(c)}_t\in \di y)\,.
$
Thus, for all $s,\delta$ non-negative and $u,w\in \R$
\begin{eqnarray*}
    \Prob(X_s^{(c)}\leq u, X_{s+\delta}^{(c)}\leq w)&=& \int_{-\infty}^u\Prob_\delta(x, (-\infty, w])\di \Prob(X_s^{(c)}\leq x)\\
    &\stackrel{s\rightarrow \infty}{\longrightarrow}&\int_{-\infty}^u\Prob_\delta(x, (-\infty, w])\di \Prob(X_\infty^{(c)}\leq x)\\
    &=&\Prob(Y_0^{(c)}\leq u, Y_{\delta}^{(c)}\leq w)\,,
\end{eqnarray*}
where in the second line we use the Portmanteau theorem because $\lim_{s\rightarrow\infty}\Prob(X_s^{(c)}\leq x)=\Prob(X_\infty^{(c)}\leq x)$
and $\Prob_\delta(x, (-\infty, w])$ is a bounded \zm{continuous} function.  By Corollary \ref{cor3} eq. (\ref{eq:stat-joint-drift}) we get the distribution
of $(Y_0^{(c)}, Y_{\delta}^{(c)})$. The proof is complete.
\qed

\subsection{Proof of Theorem \ref{thasimpBMS}}\label{thasimpBMS:proof}
Assume that $u>0$, $u>x_R$ and \zm{put $X_\infty\coloneqq X^{(0)}_\infty$}. The constant $C$ may change from line to line below.
Recall that 
$$
\Prob(X_{\infty}>x)=\left\{
\begin{array}{ll}
 1-\frac{1}{2}e^{\alpha(x-x_R)}    & \mbox{if}\,\,\,x<x_R  \\
 \frac{1}{2}e^{-\alpha(x-x_R)}    &\mbox{if} \,\,\,x\geq x_R
\end{array}
\right.
$$
for $x\in \RL$ and the epoch $S$ of the first jump of the Poisson process $N$ under condition $N(T)\geq 1$ has the density function 
\begin{equation}\label{deng}
g_T(x)=\frac{\lambda e^{-\lambda x}}{1-e^{-\lambda T}}\,\ind_{(0,T)}(x)\,.
\end{equation} 
Moreover, for $t>0$ the random variable $M_t=\sup_{s\in [0,t]}W_s$ has the density function
$f_t(x)=\frac{2}{\sqrt{2\pi t}} e^{-\frac{x^2}{2 t}}\,\ind_{(0,\infty)}(x)$.
\\
First, notice that
\begin{equation}\label{splitprob}
    \Prob(\sup_{t\in [0,T]}Y_t>u) =
    e^{-\lambda T}\Prob(\sup_{t\in [0,T]}Y_t>u | N_T=0)+(1-e^{-\lambda T})\Prob(\sup_{t\in [0,T]}Y_t>u | N_T\geq 1)\,.
\end{equation}
Let us show for the first term of (\ref{splitprob}) that
 $$
e^{-\lambda T}\Prob(\sup_{t\in [0, T]}Y_t>u | N_T=0)\sim e^{\alpha x_R}\Phi(\alpha\sigma\sqrt{T})e^{-\alpha u}\,.
 $$
Thus, we obtain
\begin{eqnarray*}
    \Prob(\sup_{t\in [0,T]}Y_t>u | N_T=0)&=& \Prob(X_\infty+\sup_{t\in [0,T]}W_t>u)\\
    &=&\int_0^{\infty}\Prob(X_\infty>u-x)f_T(x)\di x\\
    &=&\int_0^{u-x_R}\Prob(X_\infty>u-x)f_T(x)\di x+\int_{u-x_R}^{\infty}\Prob(X_\infty>u-x)f_T(x)\di x\\
    &=:& A_1+A_2\,.
\end{eqnarray*}
Let us first estimate $A_2$
\begin{eqnarray*}
    A_2&=& \int_{u-x_R}^{\infty}(1-\frac{1}{2}e^{\alpha(u-x-x_R)})f_T(x)\di x\\
    &\leq & \int_{u-x_R}^{\infty}f_T(x)\di x\\
    &=&2\Prob(W_T>u-x_R)\\
    &=&O(e^{-Cu^2})
\end{eqnarray*}
as $u\rightarrow\infty$.
Now, let us compute the term $A_1$
\begin{eqnarray*}
    A_1&=&\frac{1}{2}\int_0^{u-x_R}e^{-\alpha(u-x-x_R)}f_T(x)\di x\\
    &=&\frac{1}{2}e^{-\alpha(u-x_R)}\frac{2}{\sqrt{2\pi T}}\int_0^{u-x_R}e^{\alpha x}e^{-\frac{x^2}{2 T}}\di x\\
    &=&e^{-\alpha(u-x_R)}\frac{1}{\sqrt{2\pi T}}\int_0^{u-x_R}e^{-\frac{1}{2}(x/\sqrt{T}-\alpha \sqrt{T})^2+\frac{1}{2}\alpha^2 T}\di x\\
    &=&e^{-\alpha(u-x_R)}e^{\lambda T}\frac{1}{\sqrt{2\pi}}\int_{-\alpha\sqrt{T}}^{\frac{u-x_R}{\sqrt{T}}-\alpha\sqrt{T}}e^{-\frac{1}{2}x^2}\di x\,.
\end{eqnarray*}
 Since 
as $u\rightarrow\infty$
 $$
\frac{1}{\sqrt{2\pi}}\int_{-\alpha\sqrt{T}}^{\frac{u-x_R}{\sqrt{T}}-\alpha\sqrt{T}}e^{-\frac{1}{2}x^2}\di x
\rightarrow \Phi(\alpha\sqrt{T}),
 $$  we get that 
 $$
 A_1\sim e^{\lambda T}e^{\alpha x_R}\Phi(\alpha\sqrt{T})e^{-\alpha u},\qquad u\rightarrow\infty.
 $$ 
 
Let us consider the second term of (\ref{splitprob}), that is
 $$
B\coloneqq\Prob(\sup_{t\in [0,T]}Y_t>u | N_T\geq 1)= \Prob(\sup_{t\in [0,S]}Y_t>u \lor \sup_{t\in [S, T]}Y_t>u| N_T\geq 1)\,,
 $$
where $S$ is the epoch of the first jump of the Poisson process $N$. Notice that
\begin{equation*}
\Prob(\sup_{t\in [0,S]}Y_t>u|N_T\geq 1)\leq B\leq \Prob(\sup_{t\in [0,S]}Y_t>u|N_T\geq 1)+\Prob(\sup_{t\in [S, T]}Y_t>u| N_T\geq 1).
\end{equation*}
Consequently, it is enough to show that
\begin{equation}\label{estBu}  \Prob(\sup_{t\in [0,S]}Y_t>u|N_T\geq 1)\sim 
\frac{\lambda e^{\alpha x_R}}{1-e^{-\lambda T}}\int_0^T\Phi(\alpha\sqrt{s})\di s e^{-\alpha u}, \quad u \to \infty,
\end{equation} 
and
\begin{equation}\label{estBU} 
\Prob(\sup_{t\in [S, T]}Y_t>u| N_T\geq 1)=o(e^{-\alpha u})\,.
\end{equation}
For all $u>0$ we have  
\begin{eqnarray*}
  \Prob(\sup_{t\in [0,S]}Y_t>u|N_T\geq 1)&=&\Prob(X_\infty+\sup_{t\in [0,S]}W_t>u)\\
   &=& \int_0^T g_T(t)\di t\int_0^\infty \Prob(X_\infty>u-x)f_t(x)\di x\\
   &=& \frac{1}{2}\int_0^T g_T(t)\di t\int_0^{u-x_R}e^{-\alpha(u-x-x_R)}f_t(x)\di x\\
   &&\,\,+\int_0^T g_T(t)\di t\int_{u-x_R}^\infty[1-\frac{1}{2}e^{\alpha(u-x-x_R)}]f_t(x)\di x\\
   &=:&B_1+B_2\,.
\end{eqnarray*}

Next, we compute $B_1$
\begin{eqnarray*}
    B_1&=& \frac{1}{2}e^{\alpha(x_R-u)}\int_0^T g_T(t)\di t\int_0^{u-x_R}e^{\alpha x}\frac{2}{\sqrt{2\pi t}}
    e^{-\frac{x^2}{2 t}}\di x\\
    &=& e^{\alpha(x_R-u)}\int_0^T g_T(t)\di t\int_0^{u-x_R}\frac{1}{\sqrt{2\pi t}}
    e^{-\frac{1}{2}(x/\sqrt{t}-\alpha \sqrt{t})^2+\frac{1}{2}\alpha^2 t}\di x\\
    &=& \frac{e^{\alpha(x_R-u)}}{1-e^{-\lambda T}}\int_0^T{\lambda\di t} \int_0^{u-x_R}\frac{1}{\sqrt{2\pi t}}
    e^{-\frac{1}{2}(x/\sqrt{t}-\alpha \sqrt{t})^2}\di x
   \\
    &=& \frac{\lambda e^{\alpha(x_R-u)}}{1-e^{-\lambda T}}\int_0^T\di t\int_{-\alpha \sqrt{t}}^{\frac{u-x_R}{\sqrt{t}}-\alpha \sqrt{t}}\frac{1}{\sqrt{2\pi}}
    e^{-\frac{1}{2}x^2}\di x\\
    &\sim & e^{-\alpha u}\frac{\lambda e^{\alpha x_R}}{1-e^{-\lambda T}}\int_0^T\Phi(\alpha\sqrt{t})\di t
\end{eqnarray*}
as $u\rightarrow\infty$.
The second summand $B_2$ is estimated as follows
\begin{eqnarray*}
    B_2&\leq & \int_0^T g_T(t)\di t\int_{u-x_R}^\infty f_t(x)\di x\\
    &=& \int_0^T \Prob(M_t>u-x_R) g_T(t)\di t\\
    &=& 2\int_0^T \Prob(W_t>u-x_R) g_T(t)\di t\\
    &\leq & 2\int_0^T \Prob(W_T>u-x_R) g_T(t)\di t\\
    &=&2\Prob(W_T>u-x_R)=O(e^{-Cu^2})
\end{eqnarray*}
as $u\rightarrow\infty$.
Combining $B_1$ with $B_2$ establishes \eqref{estBu}.\\
Next,  note first that
$$
\Prob(\sup_{t\in [S, T]}Y_t>u| N_T\geq 1)=\frac{\Prob(\sup_{t\in [S, T]}Y_t>u, N_T\geq 1)}{\Prob(N_T\geq 1)}
$$
and
$$
 \Prob(\sup_{t\in [S, T]}Y_t>u, N_T\geq 1)=\sum_{k=1}^\infty \Prob(\sup_{t\in [S, T]}Y_t>u, N_T=k)\,.
$$
Let $S=U^{(1)}, \ldots, U^{(k)}$, $k\geq 1$ be the epochs of the consecutive jumps of the Poisson process $N$. 
Then, we obtain
\begin{eqnarray*}
  \lefteqn{\Prob(\sup_{t\in [S, T]}Y_t>u, N_T=k)}\\
  &=&  \Prob(\{\{\sup_{t\in [U^{(1)}, U^{(2)}]}Y_t>u\} \lor \{\sup_{t\in [U^{(2)}, U^{(3)}]}Y_t>u \}\lor\ldots \lor \{ \sup_{t\in [U^{(k)}, T]}Y_t>u\}\}\land \{N_T=k\}) \\
  &=& \Prob(\{\{x_R+\sup_{t\in [U^{(1)}, U^{(2)}]}W^{(1)}_t>u\}\lor\ldots \lor \{ x_R+\sup_{t\in [U^{(k)}, T]}W^{(k)}_t>u\}\}\land \{N_T=k\})\\
  &\leq& \Prob(\sup_{t\in [U^{(1)}, U^{(2)}]}W_t>u-x_R, N_T=k)+\ldots +\Prob(\sup_{t\in [U^{(k)}, T]}W_t>u-x_R, N_T=k)\\
  &\leq & k\Prob(\sup_{t\in [0,T]}W_t>u-x_R)\\
  &=&2k\Prob(W_T>u-x_R)=k O(e^{-Cu^2}),
\end{eqnarray*}
as $u\rightarrow\infty$, where $W^{(1)},\ldots,W^{(k)}$ are independent copies of the Wiener process $W$. So we get
\begin{equation}\label{asS1}
\Prob(\sup_{t\in [S, T]}Y_t>u, N_T\geq 1)=\sum_{k=1}^\infty e^{-\lambda\zm{T}} \frac{(\lambda \zm{T})^k}{k!}k O(e^{-Cu^2})
\end{equation}
which is $O(e^{-Cu^2})=o(e^{-\alpha u})$. Hence  \eqref{estBU} follows establishing the proof.  
\qed

\subsection{Proof of Theorem \ref{thasympSBMl}}\label{thasympSBMl:proof}
    Notice first that
    \begin{equation}\label{decompr}
    \Prob(\sup_{t\in [0,T]}{Y_t>u}, Y_T>uz)=\Prob(\sup_{t\in [0,T]}{Y_t>u})-
    \Prob(\sup_{t\in [0,T]}{Y_t>u}, Y_T\leq uz)\,.   
    \end{equation}
For the second term, we have  
\begin{eqnarray}
    \lefteqn{\Prob(\sup_{t\in [0,T]}{Y_t>u}, Y_T\leq uz)=}\label{secondt}\\
    &&e^{-\lambda T}\Prob(\sup_{t\in [0,T]}{Y_t>u}, Y_T\leq uz|N_T=0)+(1-e^{-\lambda T})\Prob(\sup_{t\in [0,T]}{Y_t>u}, Y_T\leq uz|N_T\geq 1)\nonumber\,.
\end{eqnarray}
Let $M_T=\sup_{t\in [0,T]}W_t$ and recall $Y_0=X_\infty$
being a Laplace distributed random variable \zm{($X_\infty\coloneqq X^{(0)}_\infty$)}. Then for $0< z<1$ we get
\begin{eqnarray*}
    \lefteqn{\Prob(\sup_{t\in [0,T]}{Y_t>u}, Y_T\leq uz|N_T=0)=}\\
    &&\Prob(X_\infty+M_T>u, X_\infty+W_T\leq uz)\\
    &=& \int_{-\infty}^{\infty}\Prob(M_T>u-y, W_T\leq uz-y)
    f_{X_\infty}(y)\di y\\
    &=& \int_{-\infty}^{u}\Prob(W_T\geq u(2-z)-y)f_{X_\infty}(y)\di y
    +\int_{u}^{\infty}\Prob(W_T\leq uz-y)f_{X_\infty}(y)\di y\\
    &\leq &\Prob(W_T\geq u(1-z))\int_{-\infty}^{u}f_{X_\infty}(y)\di y
    +\Prob(W_T\leq uz-u)\int_{u}^{\infty}f_{X_\infty}(y)\di y\\
 &=& \Prob(W_T\geq u(1-z))\Prob(X_\infty\leq u)+\Prob(W_T\leq u(z-1))\Prob(X_\infty>u)\\ 
&=& \Prob(W_T\ge u(1-z))=O(e^{-Cu^2})\,,
    \end{eqnarray*}
where in the third equality, we used $\Prob(M_T>m, W_T\leq y)=P(W_T\leq y)$ if $m\leq 0$ and the reflection principle that is 
$\Prob(M_T>m, W_T\leq y)=P(W_T\geq 2m-y)$ if $m>y$ and $m>0$.

For $z<0$ and $uz-x_R<0$ consider the term
\begin{eqnarray*}
    \Prob(\sup_{t\in [0,T]}{Y_t>u}, Y_T\leq uz|N_T\geq 1)&\leq&
    \Prob(Y_T\leq uz|N_T\geq 1)\\
    &=&\Prob(x_R+W_{S}\leq uz)\\
    &\leq &\Prob(W_{T}\leq uz-x_R)\\
    &=&\Prob(W_{T}\geq -uz+x_R)=O(e^{-Cu^2}),
\end{eqnarray*}
as $u\rightarrow\infty$, where $S$ is now the time since the last resetting at time $T$.  
Hence, using (\ref{decompr}) and Theorem  \ref{thasimpBMS} we obtain the asymptotic
for $z<0$.

The second term on the left side in (\ref{secondt})
is left to estimate for $0<z<1$ because the first term is $O(e^{-Cu^2})$ as shown above for all $z<1$. 
We estimate it in the following way
\begin{equation}\label{sup1}
\Prob(\sup_{t\in [0,T]}Y_t>u, Y_T\leq uz|N_T\geq 1)\geq \Prob(\sup_{t\in [0,S]}Y_t>u, Y_T\leq uz|N_T\geq 1)
\end{equation}
and
\begin{eqnarray}
\lefteqn{\Prob(\sup_{t\in [0,T]}Y_t>u, Y_T\leq uz|N_T\geq 1)\leq}\label{sup2}\\
&&\Prob(\sup_{t\in [0,S]}Y_t>u, Y_T\leq uz|N_T\geq 1)+
\Prob(\sup_{t\in [S,T]}Y_t>u, Y_T\leq uz|N_T\geq 1)\,,\nonumber
\end{eqnarray}
where $S$ is the epoch of the first resetting, which under condition $N_T\geq 1$ has the density function $g_T$ given in (\ref{deng}). Since
\begin{equation}\label{sup3}
\Prob(\sup_{t\in [S,T]}Y_t>u, Y_T\leq uz|N_T\geq 1)\leq \Prob(\sup_{t\in [S,T]}Y_t>u | N_T\geq 1)=O(e^{-Cu^2}),
\end{equation}
as $u\rightarrow\infty$
(see the derivation of (\ref{asS1})),
it is enough to consider $\Prob(\sup_{t\in [0,S]}Y_t>u, Y_T\leq uz|N_T\geq 1)$. 
For $0<s\leq T$  
$$
\Prob(X_\infty+M_s>u)=e^{\lambda s}e^{\alpha x_R}\Phi(\alpha \sqrt{s})
\exp(-\alpha u)(1+o(1)),
$$
as $u\rightarrow\infty$ where $o(1)$ can be bounded by $o(1)$ independent of $s$ (see the estimation of $A_1$ and $A_2$ in the 
proof of Theorem \ref{thasimpBMS}).
Let $S_1$ be the time since the last resetting at time $T$ then we have
\begin{eqnarray*}
    \Prob(\sup_{t\in [0,S]}Y_t>u, Y_T\leq uz|N_T\geq 1)&=& \Prob(X_\infty+M_S>u,
    W_{S_1}+x_R\leq uz)\\
    &\leq &\Prob(X_\infty+M_S>u)\\
    &=&\int_0^T\Prob(X_\infty+M_s>u)g_T(s)\di s \\
    &=&\frac{\lambda e^{\alpha x_R}e^{-\alpha u}}{1-e^{-\lambda T}}\int_0^T\Phi(\alpha\sqrt{s})\di s (1+o(1)),
\end{eqnarray*}
as $u\rightarrow\infty$.
Moreover, for $uz-x_R>0$
\begin{eqnarray*}
  \Prob(X_\infty+M_S>u, W_{S_1}+x_R\leq uz )&=&
  \int_{[0,T]^2}\Prob(X_\infty+M_s>u, W_{t}+x_R\leq uz)F(\di s,\di t)\\
  &=&\int_{[0,T]^2}\Prob(X_\infty+M_s>u)\Prob( W_{t}+x_R\leq uz)F(\di s,\di t)\\
  &\geq &\int_{[0,T]^2}\Prob(X_\infty+M_s>u)\Prob( W_{T}\leq uz-x_R)F(\di s,\di t)\\
  &= &\Prob(X_\infty+M_S>u)\Prob(W_{T}\leq uz-x_R)\\
  &=&\frac{\lambda e^{\alpha x_R}e^{-\alpha u}}{1-e^{-\lambda T}}(1+o(1))\int_0^T\Phi(\alpha\sqrt{s})\di s
\end{eqnarray*}
as $u\rightarrow\infty$, where $F(s,t)$ is the distribution function of $(S,S_1)$.
In the first equality, we use the independence of $(S,S_1)$ and $X_\infty$. In the second equality, 
we apply independence after resetting which together with the calculations above, gives that
$$
\Prob(\sup_{t\in [0,S]}Y_t>u, Y_T\leq uz|N_T\geq 1)\sim
\exp(-\alpha u)\frac{\lambda e^{\alpha x_R}}{1-e^{-\lambda T}}\int_0^T\Phi(\alpha\sqrt{s})\di s, \quad u\rightarrow\infty,
$$
which by (\ref{secondt}), (\ref{sup1}), (\ref{sup2}), (\ref{sup3})    
for $0<z<1$ gives that
$$
\Prob(\sup_{t\in [0,T]}Y_t>u, Y_T\leq uz)\sim
\exp(-\alpha u)\lambda e^{\alpha x_R}\int_0^T\Phi(\alpha\sqrt{s})\di s
$$
as $u\rightarrow\infty$. Thus, using (\ref{decompr}) and Theorem \ref{thasimpBMS} we get the asymptotic for $0<z<1$. The proof is complete.
\qed

\section*{Acknowledgments}
\noindent
K. D\c{e}bicki was partially supported by the National Science Centre, Poland,  
\\
Grant No 2024/55/B/ST1/01062
(2025-2028).   

\bibliography{Ref.bib}

@article{DebickiHashorvaNovikov2026,
  author  = {Dębicki, Krzysztof and Hashorva, Enkelejd and Novikov, Svyatoslav},
  title   = {Expected Infimum and Persistence Probabilities of Log-Normal Stationary {Brown--Resnick} Processes},
  journal = {Manuscript submitted to Extremes},
  year    = {2026}, 
}

@article{adler2014existence,
  author  = {Adler, Robert J. and Moldavskaya, Elina and Samorodnitsky, Gennady},
  title   = {On the existence of paths between points in high level excursion sets of {G}aussian random fields},
  journal = {The Annals of Probability},
  year    = {2014},
  volume  = {42},
  number  = {3},
  pages   = {1020--1053},
  doi     = {10.1214/12-AOP794}
}

@article{chakrabarty2018asymptotic,
  title     = {Asymptotic behaviour of high {G}aussian minima},
  author    = {Chakrabarty, Arijit and Samorodnitsky, Gennady},
  journal   = {Stochastic Processes and their Applications},
  volume    = {128},
  number    = {7},
  pages     = {2297--2324},
  year      = {2018},
  doi       = {10.1016/j.spa.2017.09.008},
  publisher = {Elsevier}
}

@article{dkebicki2016parisian,
  title     = {Parisian ruin over a finite-time horizon},
  author    = {Dębicki, Krzysztof and Hashorva, Enkelejd and Ji, Lanpeng},
  journal   = {Science China Mathematics},
  volume    = {59},
  number    = {3},
  pages     = {557--572},
  year      = {2016},
  doi       = {10.1007/s11425-015-5073-6},
  publisher = {Springer}
}

@article{Evans2011DiffusionResetting,
  title   = {Diffusion with Stochastic Resetting},
  author  = {M.R. Evans and S.N. Majumdar},
  journal = {Physical Review Letters},
  volume  = {106},
  number  = {16},
  pages   = {160601},
  year    = {2011},
  publisher = {American Physical Society},
  doi     = {10.1103/PhysRevLett.106.160601}
}

@article{EvansMajumdar2011JPA,
  author    = {M.R. Evans and S.N. Majumdar},
  title     = {Diffusion with optimal resetting},
  journal   = {Journal of Physics A: Mathematical and Theoretical},
  volume    = {44},
  number    = {43},
  pages     = {435001},
  year      = {2011},
  doi       = {10.1088/1751-8113/44/43/435001}
}

@incollection{Hartmann2024,
  author    = {Hartmann, Alexander K. and Majumdar, Satya N. and Schehr, Gr{\'e}gory},
  title     = {The Distribution of the Maximum of Independent Resetting {B}rownian Motions},
  booktitle = {Target Search Problems},
  editor    = {Grebenkov, Denis and Metzler, Ralf and Oshanin, Gleb},
  series    = {Modeling and Simulation in Science, Engineering and Technology},
  publisher = {Springer},
  address   = {Cham},
  year      = {2024},
  pages     = {415--435},
  isbn      = {978-3-031-67802-8},
  doi       = {10.1007/978-3-031-67802-8_15},
  url       = {https://doi.org/10.1007/978-3-031-67802-8_15}
}

@book{Viswanathan2011,
  author    = {Viswanathan, G. M. and da Luz, M. G. E. and Raposo, E. P. and Stanley, H. E.},
  title     = {The Physics of Foraging: An Introduction to Random Searches and Biological Encounters},
  year      = {2011},
  publisher = {Cambridge University Press},
  address   = {Cambridge}
}

@book{johnson2000continuous,
  title     = {Continuous Multivariate Distributions, Volume 1: Models and Applications},
  author    = {Johnson, Norman L. and Kotz, Samuel and Balakrishnan, Narayanaswamy},
  year      = {2000},
  edition   = {2nd},
  publisher = {Wiley},
  address   = {New York},
  series    = {Wiley Series in Probability and Statistics},
  isbn      = {978-0471183877}
}

@article{MagdziarzTazbierski2022,
  author       = {Magdziarz, Marcin and Taźbierski, Kacper},
  title        = {Stochastic representation of processes with resetting},
  journal      = {Physical Review E},
  volume       = {106},
  number       = {1},
  pages        = {014147},
  year         = {2022},
  doi          = {10.1103/PhysRevE.106.014147},
  publisher    = {American Physical Society}
}

@article{TazbierskiMetzlerMagdziarz2025,
  author    = {Ta\'{z}bierski, Kacper and Metzler, Ralf and Magdziarz, Marcin},
  title     = {Series representation approach to stochastic processes with complete and incomplete renewal resetting},
  journal   = {New Journal of Physics},
  volume    = {27},
  number    = {7},
  pages     = {074603},
  year      = {2025},
  doi       = {10.1088/1367-2630/adee40},
  publisher = {IOP Publishing}
}

@article{Lenzi2022Transient,
  author       = {Lenzi, M. K. and Lenzi, E. K. and Guilherme, L. M. S. and Evangelista, L. R. and Ribeiro, H. V.},
  title        = {Transient anomalous diffusion in heterogeneous media with stochastic resetting},
  journal      = {Physica A: Statistical Mechanics and its Applications},
  year         = {2022},
  volume       = {588},
  pages        = {126560},
  doi          = {10.1016/j.physa.2021.126560}
}

@article{ToledoMarin2023First,
  author       = {Toledo-Mar{\'\i}n, J. Quetzalc{\'o}atl and Boyer, Denis},
  title        = {First passage time and information of a one-dimensional {B}rownian particle with stochastic resetting to random positions},
  journal      = {Physica A: Statistical Mechanics and its Applications},
  year         = {2023},
  volume       = {625},
  pages        = {129027},
  doi          = {10.1016/j.physa.2023.129027}
}

@article{AvrachenkovPiunovskiyZhang2013,
  author  = {Avrachenkov, Konstantin E. and Piunovskiy, Alexey B. and Zhang, Yi},
  title   = {Markov Processes with Restart},
  journal = {Journal of Applied Probability},
  volume  = {50},
  number  = {4},
  pages   = {960--968},
  year    = {2013},
  doi     = {10.1239/jap/1389370093}
}

@article{EvansMajumdarSchehr2020,
  author    = {Martin R. Evans and Satya N. Majumdar and Gr{\'e}gory Schehr},
  title     = {Stochastic Resetting and Applications},
  journal   = {Journal of Physics A: Mathematical and Theoretical},
  volume    = {53},
  number    = {19},
  pages     = {193001},
  year      = {2020},
  doi       = {10.1088/1751-8121/ab7cfe},
}

@article{MonteroMasoPuigdellosasVillarroel2017,
  author    = {Miquel Montero and Axel Mas{\'o}-Puigdellosas and Javier Villarroel},
  title     = {Continuous-time random walks with reset events: Historical background and new perspectives},
  journal   = {The European Physical Journal B},
  volume    = {90},
  pages     = {176},
  year      = {2017},
  doi       = {10.1140/epjb/e2017-80348-4},
}

@article{Majumdar2020EVS,
  author    = {Majumdar, S. N. and Pal, A. and Schehr, G.},
  title     = {Extreme Value Statistics of Correlated Random Variables: A Pedagogical Review},
  journal   = {Physics Reports},
  volume    = {840},
  number    = {1},
  pages     = {1--32},
  year      = {2020},
  publisher = {Elsevier},
  doi       = {10.1016/j.physrep.2019.10.005}
}

@article{SinghPal2021ResetEVS,
  author  = {Singh, P. and Pal, A.},
  title   = {Extremal statistics for stochastic resetting systems},
  journal = {Physical Review E},
  volume  = {103},
  number  = {5},
  pages   = {052119},
  year    = {2021},
  doi     = {10.1103/PhysRevE.103.052119}
}

@article{GuptaPlataPal2020JarzynskiReset,
  author  = {Gupta, D. and Plata, C. A. and Pal, A.},
  title   = {Work Fluctuations and {J}arzynski Equality in Stochastic Resetting},
  journal = {Physical Review Letters},
  volume   = {124},
  number   = {11},
  pages    = {110608},
  year     = {2020},
  doi      = {10.1103/PhysRevLett.124.110608}
}

@article{ZhangHadjiliadis2012DrawdownsCrash,
  author    = {Zhang, H. and Hadjiliadis, O.},
  title     = {Drawdowns and the Speed of Market Crash},
  journal   = {Methodology and Computing in Applied Probability},
  year      = {2012},
  volume    = {14},
  number    = {3},
  pages     = {739--752},
  doi       = {10.1007/s11009-011-9262-7},
  publisher = {Springer},
  url       = {https://doi.org/10.1007/s11009-011-9262-7}
}

@article{Silvestrov2021DampingMC,
  author  = {Silvestrov, Dmitrii and Silvestrov, Sergei and
             Abola, Brian and Biganda, Peter S. and
             Engstr{\"o}m, Christian and Mango, Joseph M. and
             Kakuba, George},
  title   = {Perturbed Markov Chains with Damping Component},
  journal = {Methodology and Computing in Applied Probability},
  year    = {2021},
  volume  = {23},
  number  = {1},
  pages   = {369--397},
  doi     = {10.1007/s11009-020-09815-9},
  url     = {https://doi.org/10.1007/s11009-020-09815-9},
  publisher = {Springer}
}
\end{document}